\documentclass[titlepage,draft,12pt]{article} 
\usepackage[usenames,dvipsnames]{xcolor}
\usepackage{pstricks,pst-plot,pst-node} 
\usepackage{amssymb,amsthm} 
\usepackage[a4paper]{geometry}

\date{}

\newcommand{\ep}{\varepsilon}
\renewcommand{\qed}{{\penalty 10000\mbox{$\quad\Box$}}}
\newcommand{\re}{\mathbb{R}}

\newcommand{\n}{\mathbb{N}}

\newcommand{\Ehat}{\widehat{E}}

\newcommand{\istogramma}[4]{
\psframe*[linecolor=#3](0,0)(#1,#2)
\psframe[linewidth=0.5\pslinewidth,linecolor=#4](0,0)(#1,#2)
\psline[linewidth=2\pslinewidth](0,#2)(#1,#2)
}

\newcommand{\bipulse}[5]{
\istogramma{#3}{#2}{#4}{#5}
\rput(#1,0){\istogramma{-#3}{#2}{#4}{#5}}
\psline[linewidth=2\pslinewidth](#3,0)(!#1 \space #3 \space sub 0)
}

\newcommand{\tripulse}[7]{
\rput(0,0){\bipulse{#1}{#4}{#5}{#6}{#7}}
\rput(#1,0){\bipulse{#2}{#4}{#5}{#6}{#7}}
\rput(!#1\space #2 add 0){\bipulse{#3}{#4}{#5}{#6}{#7}}
}

\newtheorem{thm}{Theorem}[section]

\newtheorem{rmk}[thm]{Remark}
\newtheorem{prop}[thm]{Proposition}

\newtheorem{lemma}[thm]{Lemma}

\title{The remarkable effectiveness of time-dependent damping terms
for second order evolution equations}

\author{Marina Ghisi\vspace{1ex}\\ 
{\normalsize Universit\`a degli Studi di Pisa} \\
{\normalsize Dipartimento di Matematica}\\ 
{\normalsize PISA (Italy)}\\
{\normalsize e-mail: \texttt{ghisi@dm.unipi.it}}
\and
Massimo Gobbino\vspace{1ex}\\ 
{\normalsize Universit\`a degli Studi di Pisa} \\
{\normalsize Dipartimento di Matematica}\\ 
{\normalsize PISA (Italy)}\\  
{\normalsize e-mail: \texttt{m.gobbino@dma.unipi.it}}
\and
Alain Haraux\vspace{1ex}\\ 
{\normalsize Universit\'{e} Pierre et Marie Curie} \\
{\normalsize Laboratoire Jacques-Louis Lions}\\ 
{\normalsize PARIS (France)}\\  
{\normalsize e-mail: \texttt{haraux@ann.jussieu.fr}}}

\begin{document}
\maketitle
\begin{abstract}

	We consider a second order linear evolution equation with a 
	dissipative term multiplied by a time-dependent coefficient. Our 
	aim is to design the coefficient in such a way that \emph{all} 
	solutions decay in time as fast as possible.
	
	We discover that constant coefficients do not achieve the goal, 
	as well as time-dependent coefficients that are too big. On the 
	contrary, pulsating coefficients which alternate big and small 
	values in a suitable way prove to be more effective.
	
	Our theory applies to ordinary differential equations, systems of 
	ordinary differential equations, and partial differential 
	equations of hyperbolic type.
	
\vspace{6ex}

\noindent{\bf Mathematics Subject Classification 2010 (MSC2010):}
35B40, 49J15, 49J20.


\vspace{6ex}


\end{abstract}

 
\section{Introduction}

In this paper we consider abstract evolution equations of the form
\begin{equation}
	u''(t)+2\delta(t)u'(t)+Au(t)=0,
	\label{eqn:PDE}
\end{equation}
with initial data
\begin{equation}
	u(0)=u_{0}\in D(A^{1/2}),
	\hspace{4em}
	u'(0)=u_{1}\in H,
	\label{data}
\end{equation}
where $H$ is a Hilbert space, $A$ is a self-adjoint linear operator 
on $H$ with dense domain $D(A)$, and 
$\delta:[0,+\infty)\to[0,+\infty)$ is a measurable function. We 
always assume that the spectrum of $A$ is a finite set, or an 
unbounded increasing sequence of positive eigenvalues.

We are interested in the decay rate of all solutions to problem 
(\ref{eqn:PDE})--(\ref{data}). In particular, we are interested in 
designing the coefficient $\delta(t)$ so that \emph{all} solutions 
decay as fast as possible as $t\to +\infty$.

\paragraph{\textmd{\textit{Motivation and related literature}}}

One of the motivations of the present work originates from control
theory.  Considering an equation of the form 
$$ u''(t) + \delta Bu'(t)+Au(t)=0, $$
where $B$ is a linear operator on $H$, it is known (see for example
\cite{MR1021188}) that uniform exponential decay of the solutions in
the energy space is independent of $\delta>0$.  At this level it is
already natural to seek for an optimal damping term in the class of
multiples of $B$ and even, more generally, of operators of the form
$\Lambda B$ where $\Lambda $ is some symmetric isomorphism from $H$ to
$H$ commuting with $B$.  Nevertheless such optimization suffers from
strong limitations, as we shall see below.\medskip

In \cite{MR1356557} an attempt was done to optimize
the decay rate by a somewhat different method.  It consists in
perturbing the conservative part, namely considering an equation of 
the form
$$u''(t) + \delta Bu'(t)+Au(t) + c u(t)=0, $$ 
where $c>0$.  When $c$ and $\delta$ tend to infinity in a certain way
the decay rate can be made arbitrarily large.  The problem here is
that the equation with $c$ large can be considered as driven by the
operator $cI$ rather than $A$: the nature of the problem is altered.
Some of the results from \cite {MR1356557} were later improved in
\cite{MR3120757}.\medskip

In \cite {MR1825863} a different strategy was used in the special case
of the string equation $$ u_{tt}-u_{xx} + \delta u_t = 0. $$

The constant dissipation $\delta u_t$ was replaced by $\delta(x)u_t$, and as
$\delta(x)$ approaches the singular potential $1/x$, it was shown
that the exponential decay rate can be made as large as prescribed. 
\medskip

In the present work we try a different approach consisting in making
$\delta = \delta(t) $ time-dependent.  This approach will turn out to
be as fruitful as the $x$-dependence and basically applicable to any
self-adjoint operator $A$.  Earlier in \cite{MR2179479} (see also
\cite{Smith}), time-dependent damping terms of intermittent type
(namely with $\delta(t)$ that vanishes on some intervals up to
infinity) were shown to produce exponential decay provided some
condition on the length of the intervals of effective damping
(together with the maxima and the minima of $\delta$ on those
intervals) is satisfied.  This result was extended in~\cite{np} to
some cases where the damping operator involves a delay term.  This was
one more motivation to examine the case of time-dependent damping.

\paragraph{\textmd{\textit{A simple toy model: ODEs with constant dissipation}}}

We begin our investigation by recalling the behavior of solutions in
the simplest example where $H=\re$ and $\delta(t)$ is constant, so
that (\ref{eqn:PDE}) reduces to the ordinary differential equation
\begin{equation}
	u''(t)+2\delta u'(t)+\lambda^{2}u(t)=0,
	\label{ODE-trivial}
\end{equation}
where $\delta$ and $\lambda$ are positive parameters.  This equation
can be explicitly integrated.  It turns out that the asymptotic
behavior of solutions depends on the real part of the roots of the
characteristic equation
\begin{equation}
	x^{2}+2\delta x+\lambda^{2}=0.
	\label{eqn:char}
\end{equation}

When $\delta<\lambda$, the characteristic equation (\ref{eqn:char}) 
has two complex conjugate roots with real part equal to $-\delta$. As 
a consequence, all nonzero solutions to (\ref{ODE-trivial}) decay as 
$e^{-\delta t}$. When $\delta>\lambda$, the characteristic equation 
(\ref{eqn:char}) has two real roots $r_{1}$ and $r_{2}$, with
$$r_{1}:=-\delta-\sqrt{\delta^{2}-\lambda^{2}}\sim -2\delta,
\hspace{5em}
r_{2}:=-\delta+\sqrt{\delta^{2}-\lambda^{2}}\sim
-\frac{\lambda^{2}}{2\delta}.$$

Every solution to (\ref{ODE-trivial}) is a linear combination of
$e^{-r_{1}t}$ and $e^{-r_{2}t}$.  Since we are looking for a decay
rate valid for \emph{all} solutions, in this case $e^{-r_{2}t}$ is the
best possible estimate, and it is also optimal for the generic
solution to (\ref{ODE-trivial}).  We recall that \emph{for linear
equations the slowest behavior is always generic}, and therefore the
best estimate valid for all solutions is sharp for the generic
solution.

A similar argument shows that for $\delta=\lambda$ all solutions to
(\ref{ODE-trivial}) are linear combinations of $e^{-\lambda t}$ and
$te^{-\lambda t}$, so that $te^{-\lambda t}$ is the estimate valid 
for all solutions in that case.

\begin{figure}[htbp]
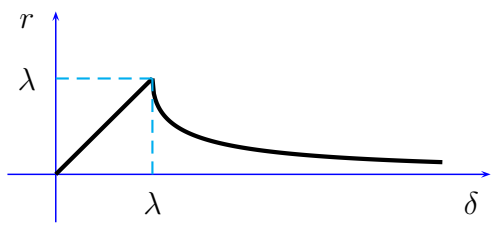

	\centering
	\psset{unit=7ex}
	\pspicture(-1,-0.3)(5,2)
	\psline[linewidth=0.7\pslinewidth, linecolor=blue]{->}(-0.5,0)(4.5,0)
	\psline[linewidth=0.7\pslinewidth, linecolor=blue]{->}(0,-0.5)(0,1.7)
	\psplot[linewidth=2\pslinewidth]{0}{1}{x}
	\psplot[plotpoints=200,linewidth=2\pslinewidth]{1}{4}{x x 2 exp 1 sub sqrt sub}
	\psline[linestyle=dashed, linecolor=cyan](1,0)(1,1)
	\psline[linestyle=dashed, linecolor=cyan](0,1)(1,1)
	\rput(4.3,-0.3){$\delta$}
	\rput(-0.3,1.6){$r$}
	\rput(1,-0.3){$\lambda$}
	\rput(-0.3,1){$\lambda$}
	\endpspicture
	\caption{Decay rate $r$ as a function of the constant dissipation $\delta$}
	\label{fig:r-delta}
\end{figure}

The graph in Figure~\ref{fig:r-delta} represents the exponent $r$ in
the optimal decay rate $e^{-rt}$ as a function of $\delta$ (namely
$r=\delta$ when $\delta<\lambda$ and $r=r_{2}$ when $\delta>\lambda$).

It clarifies that, if we limit ourselves to constant coefficients, we
cannot hope that all solutions to (\ref{ODE-trivial}) decay better
than $e^{-\lambda t}$, and actually neither better than $te^{-\lambda
t}$.  Moreover, beyond the threshold $\lambda$, the larger is $\delta$
the worse is the decay rate.  This shows also that in the general
setting of equation (\ref{eqn:PDE}), if we restrict ourselves to
constant damping coefficients $\delta(t)$, the best possible decay
rate valid for all solutions is $te^{-\nu t}$, where $\nu^{2}$ is the
smallest eigenvalue of the operator $A$.  In~\cite{MR1825863} this
remark was summed up effectively by saying that ``more is not
better''.

\paragraph{\textmd{\textit{ODEs with nonconstant dissipation}}}

The first nontrivial case we consider is the ordinary differential 
equation
\begin{equation}
	u''(t)+2\delta(t)u'(t)+\lambda^{2}u(t)=0,
	\label{eqn:ODE}
\end{equation}
where now the coefficient $\delta(t)$ is time-dependent. As long as 
$\delta\in L^{1}_{\mbox{{\scriptsize loc}}}((0,+\infty))$ the energy 
of a non-trivial solution can not vanish at any finite time. 
Moreover, if we limit
ourselves to coefficients $\delta(t)\geq\lambda$, once again there
exists a solution which decays at most as $te^{-\lambda t}$, exactly
as in the case where $\delta(t)$ is constant (see
Proposition~\ref{prop:ODE-slow}).  In other words, once again ``more
is not better'' and the overdamping prevents a faster stabilization.

Things change when we consider damping coefficients $\delta(t)$ 
which alternate intervals where they are big and intervals where they 
are small (below $\lambda$). We obtain two results.
\begin{itemize}
	\item  
	In Theorem~\ref{thm:ODE-fixed} we prove that every exponential
	decay rate can be achieved through a periodic damping coefficient.
	More precisely, for every real number $R$ there exists a periodic
	function $\delta(t)$ for which all solutions to (\ref{eqn:ODE})
	decay at least as $e^{-Rt}$.  A possible choice for the period of
	$\delta(t)$ is
	\begin{equation}
		t_{0}=\frac{\pi}{2\lambda},
		\label{defn:t0-ODE}
	\end{equation}
	hence it does \emph{not} depend on $R$.  We can also ask further
	requirements on $\delta(t)$, for example being of class $C^{\infty}$,
	or taking alternatively only two values, 0 and a sufficiently large
	positive number $K$.

	\item  
	In Theorem~\ref{thm:ODE-any} we obtain even better decay rates.
	Indeed, we prove that for every nonincreasing function
	$\varphi:[0,+\infty)\to(0,+\infty)$ one can design $\delta(t)$ in
	such a way that all solutions to (\ref{eqn:ODE}) decay at least as
	$\varphi(t)$.  In this case $\delta(t)$ is necessarily
	non-periodic and unbounded, but one can choose it of class
	$C^{\infty}$ or piecewise constant.  The proof of this second
	result relies on the first one, and the key point is that in the
	first result we can achieve any exponential decay rate $e^{-Rt}$
	with a coefficient whose period does not depend on $R$.
\end{itemize}

The results for the single equation can be easily extended to systems 
of the form 
\begin{equation}
	u_{k}''(t)+2\delta(t)u_{k}'(t)+\lambda_{k}^{2}u_{k}(t)=0
	\quad\quad\quad
	k=1,\ldots,n.
	\label{eqn:ODE-system}
\end{equation}

The only difference is that in the first result (the one with a fixed 
exponential decay rate $e^{-Rt}$) the period of $\delta(t)$ is now 
\begin{equation}
	t_{0}=\frac{\pi}{2}
	\sum_{k=1}^{n}\frac{1}{\lambda_{k}}.
	\label{defn:t0-system}
\end{equation}

Once again the period is independent of $R$, and this is the key point
to achieve any given decay rate through a non-periodic and unbounded
coefficient.  We refer to Theorem~\ref{thm:system-fixed} and
Theorem~\ref{thm:system-any} for the details.

The results obtained for (systems of) ordinary differential equations 
can be extended word for word to the Hilbert setting of 
(\ref{eqn:PDE}), provided that $A$ has a finite number of 
eigenvalues, even with infinite dimensional eigenspaces.

\paragraph{\textmd{\textit{PDEs with nonconstant dissipation}}}

It remains to consider the case of operators with an infinite number
of eigenvalues.  Due to our assumption on the spectrum, in this case
$H$ admits an orthonormal system made by eigenvectors of $A$, and
therefore the evolution equation (\ref{eqn:PDE}) is equivalent to a
system of countably many ordinary differential equations.  Looking
at~(\ref{defn:t0-system}) one could naturally guess that our theory
extends to the general setting when the series of $1/\lambda_{k}$ is
convergent.  This is actually true, but not so interesting because in
most applications the series is divergent (with the notable exception
of the beam equation, see section~\ref{sec:applications}).  This leads us
to follow a partially different path.

In Theorem~\ref{thm:PDE-fixed} we prove once again that any 
exponential decay rate $e^{-Rt}$ can be achieved through a periodic 
damping coefficient. The main difference is that now the period of 
$\delta(t)$ is 
\begin{equation}
	t_{0}=\pi\sum_{\lambda_{k}^{2}\leq 2(R+\lambda_{1})^{2}}
	\frac{1}{\lambda_{k}},
	\label{defn:t0-PDE}
\end{equation}
and hence it does depend on $R$. In analogy with the previous 
results, we can again ask further structure on $\delta(t)$, for 
example being of class $C^{\infty}$, or taking alternatively only 
three values (instead of two).

The fact that the period depends on $R$ complicates the search for
better decay rates.  Our best result is stated in
Theorem~\ref{thm:PDE-top}, where we prove that there exists
$\delta(t)$ such that all solutions to (\ref{eqn:PDE}) decay at least
as a nonincreasing function $\varphi:[0,+\infty)\to(0,1)$ that tends
to zero faster than all exponentials.  This universal decay rate is
independent of the solution, but it does depend on the operator $A$,
more precisely on its spectrum.

\paragraph{\textmd{\textit{More requirements on the damping 
coefficient}}}

The coefficients introduced in the proofs of the results quoted so far
alternate intervals where they are close to 0 and intervals where they
are very large.  Even the coefficients of class $C^{\infty}$ are just
smooth approximations of the discrete ones.  On the other hand, we
already know that both large values and values below $\lambda$ are
needed if we want fast decay rates (see Proposition~\ref{prop:ODE-slow}).

In the last part of the paper we ask ourselves whether it is essential
that the coefficient approaches 0 or exhibits sudden oscillations
between big and small values.  The answer to both questions is
negative.  In the case of the ordinary differential equation
(\ref{eqn:ODE}) we show that we can achieve any exponential decay rate
$e^{-Rt}$ through a periodic coefficient $\delta(t)$ which is always
greater than or equal to $\lambda-\ep$ and has Lipschitz constant
equal to $\ep$ (where $\ep$ is a fixed parameter).  Of course the
period of the coefficient now depends on $\ep$ and $R$. We refer to 
Theorem~\ref{thm:ODE-fixed-lip} for the details.

\paragraph{\textmd{\textit{Perspectives and open problems}}}

We consider this paper as a starting point of a research project.
Several related questions are not addressed here but could probably
deserve future investigations.  Just to give some examples, we mention
finding the optimal decay rates that can be achieved through damping
coefficients with reasonable restrictions, proving or disproving that
in the infinite dimensional setting there is a bound on the decay rate
one can achieve, extending if possible some parts of the theory to
operators with continuum spectrum, proving or disproving that random
damping coefficients are ineffective and never better than constant
ones.

\paragraph{\textmd{\textit{Structure of the paper}}}

This paper is organized as follows. In section~\ref{sec:statements} 
we state all our results. In section~\ref{sec:heuristics} we  
give a rough explanation of why ``pulsating is better''. In 
section~\ref{sec:proofs} we provide rigorous proofs. In 
section~\ref{sec:applications} we present some simple applications to 
partial differential equations.

\setcounter{equation}{0}
\section{Statements}\label{sec:statements}

For the sake of clarity we present our results in increasing order of
complexity.  We start with ordinary differential equations, we
continue with systems of ordinary differential equations, and finally
we consider the more general Hilbert setting.  In the last subsection
we investigate the same problems with additional constraints on the
damping coefficients. 

In the sequel $(t-t_{0})^{+}$ stands for 
$\max\{t-t_{0},0\}$.

\subsection{Ordinary differential equations}\label{sec:ODE}

To begin with, we consider the ordinary differential equation 
(\ref{eqn:ODE}).
In the first result we achieve any given exponential decay through a 
periodic damping coefficient. 

\begin{thm}[Single ODE, fixed exponential decay rate]\label{thm:ODE-fixed}
	Let $\lambda$ and $R$ be positive real numbers, and let $t_{0}$ 
	be defined by (\ref{defn:t0-ODE}). 
	
	Then there exists a $t_{0}$-periodic damping coefficient
	$\delta:[0,+\infty)\to[0,+\infty)$ (which one can choose either of class 
	$C^{\infty}$ or piecewise constant) such that every solution $u(t)$
	to~(\ref{eqn:ODE}) satisfies
	\begin{equation}
		|u'(t)|^{2}+\lambda^{2}|u(t)|^{2}\leq
		\left(|u'(0)|^{2}+\lambda^{2}|u(0)|^{2}\right)
		\exp\left(-R(t-t_{0})^{+}\right)
		\quad\quad
		\forall t\geq 0.
		\label{th:ODE-fixed}
	\end{equation}
	
\end{thm}

\begin{figure}[thbp]
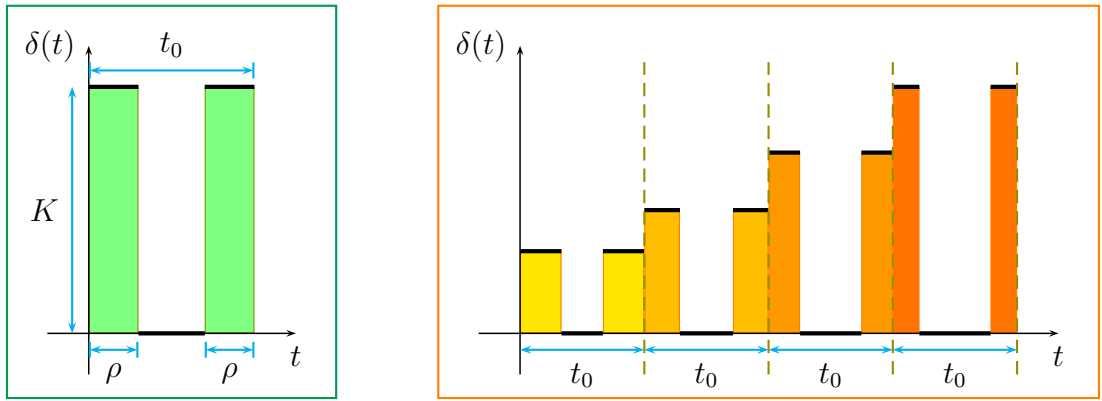

	\centering
	\psset{unit=6ex}
	\SpecialCoor
	\hfill
	\pspicture(-1,-0.8)(3,4)
	\psframe[linecolor=ForestGreen](-1,-0.8)(3,4)
	\newrgbcolor{verdino}{0.5 1 0.5}
	\bipulse{2}{3}{0.6}{verdino}{olive}
	\psline[linewidth=0.7\pslinewidth]{->}(-0.5,0)(2.5,0)
	\psline[linewidth=0.7\pslinewidth]{->}(0,-0.5)(0,3.5)
	\uput[l](0,3.5){$\delta(t)$}
	\uput[d](2.5,0){$t$}
	\pcline[offset=-0.2,linecolor=cyan]{|<->|}(0,0)(0.6,0)
	\lput{U}{\uput[270](0,0){$\rho$}}
	\pcline[offset=-0.2,linecolor=cyan]{|<->|}(1.4,0)(2,0)
	\lput{U}{\uput[270](0,0){$\rho$}}
	\pcline[offset=0.2,linecolor=cyan]{|<->|}(0,3)(2,3)
	\lput{U}{\uput[90](0,0){$t_{0}$}}
	\pcline[offset=0.2,linecolor=cyan]{<->}(0,0)(0,3)
	\lput{U}{\uput[180](0,0){$K$}}
	\endpspicture
	\hfill\hfill
	\pspicture(-1,-0.8)(7,4)
	\psframe[linecolor=orange](-1,-0.8)(7,4)
	\newrgbcolor{orange1}{1 0.9 0}
	\newrgbcolor{orange2}{1 0.75 0}
	\newrgbcolor{orange3}{1 0.6 0}
	\newrgbcolor{orange4}{1 0.45 0}
	\rput(0,0){\bipulse{1.5}{1}{0.5}{orange1}{orange}}
	\rput(1.5,0){\bipulse{1.5}{1.5}{0.43}{orange2}{orange}}
	\rput(3,0){\bipulse{1.5}{2.2}{0.38}{orange3}{orange}}
	\rput(4.5,0){\bipulse{1.5}{3}{0.32}{orange4}{orange}}
	\psline[linewidth=0.7\pslinewidth]{->}(-0.5,0)(6.5,0)
	\psline[linewidth=0.7\pslinewidth]{->}(0,-0.5)(0,3.5)
	\uput[l](0,3.5){$\delta(t)$}
	\uput[d](6.5,0){$t$}
	\multirput(0,0)(1.5,0){4}{
	\psline[linestyle=dashed,linecolor=olive](1.5,-0.5)(1.5,3.3)
	\pcline[offset=-0.2,linecolor=cyan]{<->}(0,0)(1.5,0)
	\lput{U}{\uput[270](0,0){$t_{0}$}}
	}
	\endpspicture
	\hfill\mbox{}
	\caption{possible profiles of $\delta(t)$ in
	Theorem~\ref{thm:ODE-fixed} (left) and Theorem~\ref{thm:ODE-any}
	(right)}
	\label{fig:delta-ODE}
\end{figure}

The typical profile of a damping coefficient realizing a fixed
exponential decay rate is shown in Figure~\ref{fig:delta-ODE} on the
left.  In each period there are two impulses of suitable height $K$
and duration $\rho$, one at the beginning and one at the end of the
period.  For the rest of the time the damping coefficient vanishes,
which means that there is no dissipation.  Of course, when the 
coefficient is extended by periodicity, the impulse at
the end of each period continues with the impulse at the beginning of
next period, thus giving rise to a single impulse with double
time-length.

These piecewise constant coefficients with only two values sound good 
for applications. In any case, since (\ref{eqn:ODE}) is stable under 
$L^{2}$ perturbations of the coefficient (see 
Lemma~\ref{lemma:approx}), the same effect can be achieved through a 
smooth approximation of the piecewise constant coefficient.

In the second result we show that every fixed decay rate can be 
achieved if we are allowed to exploit non-periodic and unbounded 
damping coefficients.

\begin{thm}[Single ODE, any given decay rate]\label{thm:ODE-any}
	Let $\lambda$ be a positive real number, let
	$\varphi:[0,+\infty)\to(0,+\infty)$ be a nonincreasing function,
	and let $t_{0}$ be defined by (\ref{defn:t0-ODE}).
	
	Then there exists a damping coefficient
	$\delta:[0,+\infty)\to[0,+\infty)$ (which one can choose either of class 
	$C^{\infty}$ or piecewise constant) such that every solution $u(t)$
	to (\ref{eqn:ODE}) satisfies
	\begin{equation}
		|u'(t)|^{2}+\lambda^{2}|u(t)|^{2}\leq
		\left(|u'(0)|^{2}+\lambda^{2}|u(0)|^{2}\right)
		\cdot\varphi(t)
		\quad\quad
		\forall t\geq t_{0}.
		\label{th:ODE-any}
	\end{equation}
	
\end{thm}

The typical profile of a damping coefficient realizing a given decay
rate in Theorem~\ref{thm:ODE-any} is shown in
Figure~\ref{fig:delta-ODE} on the right.  If consists in a sequence of
blocks of the type described after Theorem~\ref{thm:ODE-fixed}, the
only difference being that now the values of the parameters $K$ and
$\rho$ are different in each block.

\subsection{Systems of ordinary differential equations}\label{sec:systems}

The results for a single ordinary differential equation 
can be extended to systems. The following statement is the 
generalization of Theorem~\ref{thm:ODE-fixed}.

\begin{thm}[System of ODEs, fixed exponential decay rate]\label{thm:system-fixed}
	Let $n$ be a positive integer, let
	$(\lambda_{1},\ldots,\lambda_{n})\in(0,+\infty)^{n}$, let $R$ be a
	positive real number, and let $t_{0}$ be defined by
	(\ref{defn:t0-system}).
	
	Then there exists a $t_{0}$-periodic damping coefficient
	$\delta:[0,+\infty)\to[0,+\infty)$ (which one can choose either of
	class $C^{\infty}$ or piecewise constant) such that every solution
	$(u_{1}(t),\ldots,u_{n}(t))$ to system (\ref{eqn:ODE-system})
	satisfies
	$$\sum_{k=1}^{n}
	\left(|u_{k}'(t)|^{2}+\lambda_{k}^{2}|u_{k}(t)|^{2}\right)\leq
	\left[\sum_{k=1}^{n}
	\left(|u_{k}'(0)|^{2}+\lambda_{k}^{2}|u_{k}(0)|^{2}\right)\right]
	\exp\left(-R(t-t_{0})^{+}\right)
	$$
	for every $t\geq 0$.

\end{thm}

\begin{figure}[htbp]
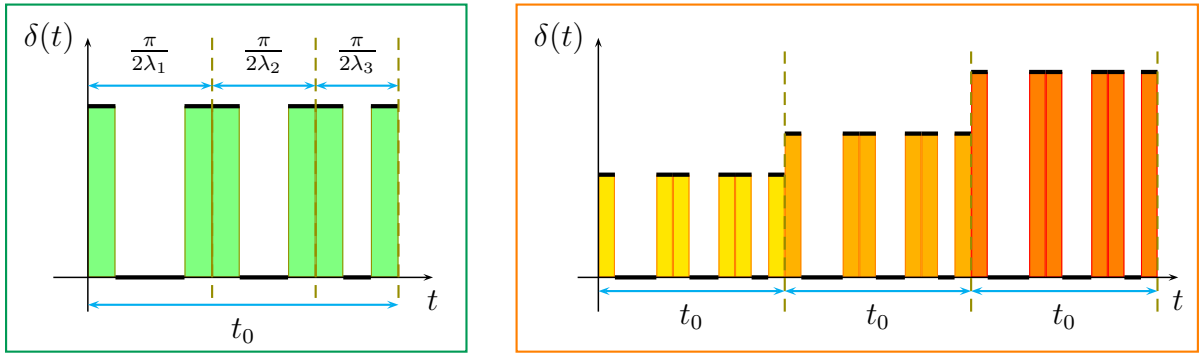

	\centering
	\psset{unit=5ex}
	\SpecialCoor
	\pspicture(-1.2,-1.1)(5.5,4)
	\psframe[linecolor=ForestGreen](-1.2,-1.1)(5.5,4)

	\tripulse{1.8}{1.5}{1.2}{2.5}{0.4}{verdino}{olive}

	\psline[linewidth=0.7\pslinewidth]{->}(-0.5,0)(5,0)
	\psline[linewidth=0.7\pslinewidth]{->}(0,-0.5)(0,3.5)
	\uput[l](0,3.5){$\delta(t)$}
	\uput[d](5,0){$t$}

	\psline[linestyle=dashed,linecolor=olive](1.8,-0.3)(1.8,3.5)
	\psline[linestyle=dashed,linecolor=olive](3.3,-0.3)(3.3,3.5)
	\psline[linestyle=dashed,linecolor=olive](4.5,-0.3)(4.5,3.5)

	\pcline[offset=-0.4,linecolor=cyan]{<->}(0,0)(4.5,0)
	\lput{U}{\uput[270](0,0){$t_{0}$}}
	\pcline[linecolor=cyan]{<->}(0,2.8)(1.8,2.8)
	\lput{U}{\uput[90](0,0){$\frac{\pi}{2\lambda_{1}}$}}
	\pcline[linecolor=cyan]{<->}(1.8,2.8)(3.3,2.8)
	\lput{U}{\uput[90](0,0){$\frac{\pi}{2\lambda_{2}}$}}
	\pcline[linecolor=cyan]{<->}(3.3,2.8)(4.5,2.8)
	\lput{U}{\uput[90](0,0){$\frac{\pi}{2\lambda_{3}}$}}

	\endpspicture
	\hfill
	\psset{xunit=3ex}
	\pspicture(-2,-1.1)(14.5,4)
	\psframe[linecolor=orange](-2,-1.1)(14.5,4)
	\newrgbcolor{orange25}{1 0.7 0}
	\newrgbcolor{orange35}{1 0.5 0}
	\rput(0,0){\tripulse{1.8}{1.5}{1.2}{1.5}{0.4}{orange1}{orange}}
	\rput(4.5,0){\tripulse{1.8}{1.5}{1.2}{2.1}{0.4}{orange25}{orange}}
	\rput(9,0){\tripulse{1.8}{1.5}{1.2}{3}{0.4}{orange35}{red}}

	\psline[linewidth=0.7\pslinewidth]{->}(-0.5,0)(14,0)
	\psline[linewidth=0.7\pslinewidth]{->}(0,-0.5)(0,3.5)
	\uput[l](0,3.5){$\delta(t)$}
	\uput[d](14,0){$t$}

	\multirput(0,0)(4.5,0){3}{
	\psline[linestyle=dashed,linecolor=olive](4.5,-0.5)(4.5,3.3)
	\pcline[offset=-0.2,linecolor=cyan]{<->}(0,0)(4.5,0)
	\lput{U}{\uput[270](0,0){$t_{0}$}}
	}

	\endpspicture
	\caption{possible profiles of $\delta(t)$ in
	Theorem~\ref{thm:system-fixed} (left) and Theorem~\ref{thm:system-any}
	(right)}
	\label{fig:delta-system}
\end{figure}

The typical profile of a damping coefficient realizing a fixed
exponential decay rate for a system is shown in
Figure~\ref{fig:delta-system} on the left.  Now the period $[0,t_{0}]$
is the union of $k$ subintervals ($k=3$ in the figure), where the
$i$-th subinterval has length $\pi/(2\lambda_{i})$.  In each
subinterval we exploit once again the profile with two impulses
described after Theorem~\ref{thm:ODE-fixed}.  We can assume that in
all subintervals the values of the parameters $K$ and $\rho$ are the
same (and hence the time-length of the vanishing phase is different),
so we end up once again with a damping coefficient with just two
values.

The following statement is the generalization of
Theorem~\ref{thm:ODE-any} to systems.

\begin{thm}[System of ODEs, any given decay rate]\label{thm:system-any}
	Let $n$ be a positive integer, let
	$(\lambda_{1},\ldots,\lambda_{n})\in(0,+\infty)^{n}$, let
	$\varphi:[0,+\infty)\to(0,+\infty)$ be a nonincreasing function,
	and let $t_{0}$ be defined by (\ref{defn:t0-system}).

	Then there exists a damping coefficient
	$\delta:[0,+\infty)\to[0,+\infty)$ (which one can choose either of
	class $C^{\infty}$ or piecewise constant) such that every solution
	$(u_{1}(t),\ldots,u_{n}(t))$ to system (\ref{eqn:ODE-system})
	satisfies
	$$\sum_{k=1}^{n}
	\left(|u_{k}'(t)|^{2}+\lambda_{k}^{2}|u_{k}(t)|^{2}\right)\leq
	\left[\sum_{k=1}^{n}
	\left(|u_{k}'(0)|^{2}+\lambda_{k}^{2}|u_{k}(0)|^{2}\right)
	\right]\cdot\varphi(t)
	\quad\quad
	\forall t\geq t_{0}.$$
	
\end{thm}

As in the case of a single equation, the typical profile of a damping
coefficient realizing a given decay rate for a system is a sequence of
blocks of the same type used in order to realize exponential decay
rates for the same system, just with different values of the
parameters $K$ and $\rho$ in different blocks (see
Figure~\ref{fig:delta-system} on the right).

\subsection{Partial differential equations}\label{sec:PDE}

We examine now equation (\ref{eqn:PDE}) in the general Hilbert 
setting. As usual, we consider weak solutions with regularity
$$u\in C^{0}\left([0,+\infty),D(A^{1/2})\right) \cap
C^{1}([0,+\infty),H).$$

In the first result we achieve once again a given exponential decay
rate through a periodic damping coefficient.  In contrast with the
case of ordinary differential equations or systems, the period of the
coefficient now does depend on the decay rate.

\begin{thm}[PDE, fixed exponential decay rate]\label{thm:PDE-fixed}
	Let $H$ be a Hilbert space, and let $A$ be a self-adjoint
	nonnegative operator on $H$ with dense domain $D(A)$.  Let us
	assume that the spectrum of $A$ is an increasing unbounded
	sequence of positive real numbers $\{\lambda_{k}^{2}\}_{k\geq 1}$
	(with the agreement that $\lambda_{k}>0$ for every $k\geq 1$).
	
	Let $R$ be a positive real number, 
	and let $t_{0}$ be defined by (\ref{defn:t0-PDE}).  
	
	Then there exists a $t_{0}$-periodic damping coefficient
	$\delta:[0,+\infty)\to[0,+\infty)$ (which one can choose either of
	class $C^{\infty}$ or piecewise constant) such that every weak solution
	to equation (\ref{eqn:PDE}) satisfies
	\begin{equation}
		|u'(t)|^{2}+|A^{1/2}u(t)|^{2}\leq
		\left(|u'(0)|^{2}+|A^{1/2}u(0)|^{2}\right)
		\exp\left(-R(t-t_{0})^{+}\right)
		\label{th:PDE-fixed}
	\end{equation}
	for every $t\geq 0$.

\end{thm}

\begin{figure}[htbp]
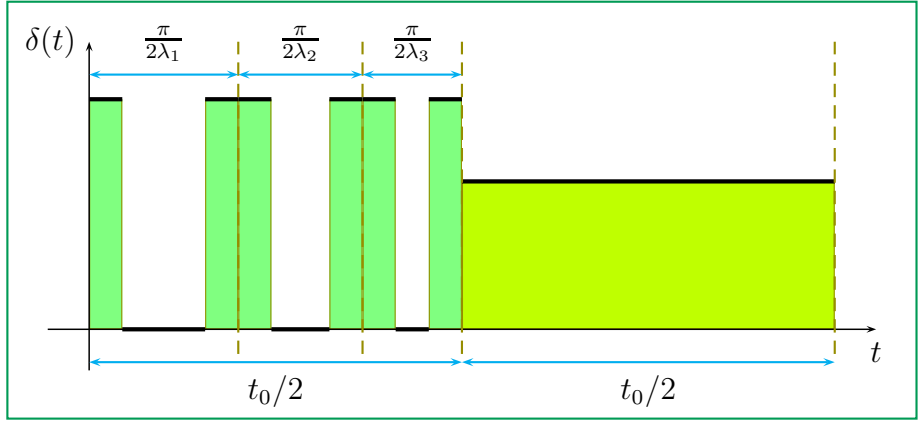

	\centering
	\SpecialCoor
	\psset{unit=6ex}
	\pspicture(-1,-1.1)(10,4)
	\psframe[linecolor=ForestGreen](-1,-1.1)(10,4)
	\newrgbcolor{verdino2}{0.75 1 0}
	\tripulse{1.8}{1.5}{1.2}{2.8}{0.4}{verdino}{olive}
	\rput(4.5,0){\istogramma{4.5}{1.8}{verdino2}{olive}}

	\psline[linewidth=0.7\pslinewidth]{->}(-0.5,0)(9.5,0)
	\psline[linewidth=0.7\pslinewidth]{->}(0,-0.5)(0,3.5)
	\uput[l](0,3.5){$\delta(t)$}
	\uput[d](9.5,0){$t$}

	\psline[linestyle=dashed,linecolor=olive](1.8,-0.3)(1.8,3.5)
	\psline[linestyle=dashed,linecolor=olive](3.3,-0.3)(3.3,3.5)
	\psline[linestyle=dashed,linecolor=olive](4.5,-0.3)(4.5,3.5)
	\psline[linestyle=dashed,linecolor=olive](9,-0.3)(9,3.5)

	\pcline[offset=-0.4,linecolor=cyan]{<->}(0,0)(4.5,0)
	\lput{U}{\uput[270](0,0){$t_{0}/2$}}
	\pcline[offset=-0.4,linecolor=cyan]{<->}(4.5,0)(9,0)
	\lput{U}{\uput[270](0,0){$t_{0}/2$}}
	\pcline[linecolor=cyan]{<->}(0,3.1)(1.8,3.1)
	\lput{U}{\uput[90](0,0){$\frac{\pi}{2\lambda_{1}}$}}
	\pcline[linecolor=cyan]{<->}(1.8,3.1)(3.3,3.1)
	\lput{U}{\uput[90](0,0){$\frac{\pi}{2\lambda_{2}}$}}
	\pcline[linecolor=cyan]{<->}(3.3,3.1)(4.5,3.1)
	\lput{U}{\uput[90](0,0){$\frac{\pi}{2\lambda_{3}}$}}

	\endpspicture
	\caption{possible profiles of $\delta(t)$ in Theorem~\ref{thm:PDE-fixed}}
	\label{fig:delta-PDE}
\end{figure}

The typical profile of a damping coefficient realizing a fixed
exponential decay rate for the full equation (\ref{eqn:PDE}) is shown
in Figure~\ref{fig:delta-PDE}.  Now the period $[0,t_{0}]$ is divided
into two subintervals of the same length $t_{0}/2$.  In the second
subinterval the damping coefficient is constant.  The first half of
the period is in turn divided into subintervals of length
$\pi/(2\lambda_{i})$, where the $\lambda_{i}$'s are those which
contribute to the sum in the right-hand side of (\ref{defn:t0-PDE}).
In each of these subintervals we have again the same profile with two
impulses described after Theorem~\ref{thm:ODE-fixed}.  We can assume
that the parameters $K$ and $\rho$ are the same in all subintervals of
the first half of the period, but the constant in the second half of
the period might differ from $K$.  Therefore, now the damping
coefficient takes in general three values instead of two.

In the second result as usual we allow ourselves to exploit
non-periodic and unbounded damping coefficients.  What we obtain is a
decay rate which is faster than all exponentials. This decay rate 
does depend on the operator.

\begin{thm}[PDE, decay rate faster than all	exponentials]\label{thm:PDE-top}
	Let the Hilbert space $H$ and the operator $A$ be as in 
	Theorem~\ref{thm:PDE-fixed}. 
	
	Then there exist a damping coefficient
	$\delta:[0,+\infty)\to[0,+\infty)$ (which one can choose either of
	class $C^{\infty}$ or piecewise constant) and a nonincreasing function
	$\varphi:[0,+\infty)\to(0,+\infty)$ such that
	\begin{equation}
		\lim_{t\to +\infty}\varphi(t)e^{Rt}=0
		\quad\quad
		\forall R>0,
		\label{th:phi-decay}
	\end{equation}
	and such that every weak solution to equation (\ref{eqn:PDE})
	satisfies
	\begin{equation}
		|u'(t)|^{2}+|A^{1/2}u(t)|^{2}\leq
		\left(|u'(0)|^{2}+|A^{1/2}u(0)|^{2}\right)
		\cdot\varphi(t)
		\quad\quad
		\forall t\geq 0.
		\label{th:PDE-top}
	\end{equation}
	
\end{thm}

As in the case of ordinary differential equations or systems, the
standard profile of a damping coefficient provided by
Theorem~\ref{thm:PDE-top} is a sequence of blocks of the same type as
those introduced for Theorem~\ref{thm:PDE-fixed}.  The main difference
is that now also the time-length $t_{0}$ is different in different
blocks, and increases with time.  The reason is that now $t_{0}$
depends on $R$, and when $R$ increases a larger number of eigenvalues
contributes to (\ref{defn:t0-PDE}), and hence $t_{0}$ increases as
well.

\subsection{Further requirements on the damping coefficient}\label{sec:lip}

In the last part of the paper we investigate which features of the
damping coefficient are essential when one wants to achieve fast decay
rates.  To begin with, in the following statement we list three simple
situations where the decay rate of solutions can be bounded \emph{from
below}.  We point out that, when looking for estimates from below,
there is almost no loss of generality in considering just the ordinary
differential equation~(\ref{eqn:ODE}).

\begin{prop}[Estimates of the decay rate from below]\label{prop:ODE-slow}
	Let us consider equation~(\ref{eqn:ODE}) for some positive real 
	number $\lambda$ and some nonnegative damping coefficient 
	$\delta\in L^{1}_{loc}((0,+\infty))$.
	\begin{enumerate}
		\renewcommand{\labelenumi}{(\arabic{enumi})}
		\item  If $\delta\in L^{1}((0,+\infty))$, then all 
		nonzero solutions do not decay to zero.
	
		\item  If there exists a constant $M>0$ such that 
		$\delta(t)\leq M$ for every $t\geq 0$, then all 
		solutions satisfy
		$$|u'(t)|^{2}+\lambda^{2}|u(t)|^{2}\geq
		\left(|u'(0)|^{2}+\lambda^{2}|u(0)|^{2}\right)
		e^{-4Mt}
		\quad\quad
		\forall t\geq 0.$$
	
		\item If there exists $T\geq 0$ such that
		$\delta(t)\geq\lambda$ for every $t\geq T$, then there
		exist $T_{*}\geq T$ and a solution $u(t)$ to
		(\ref{eqn:ODE}) such that
		\begin{equation}
			|u(t)|\geq te^{-\lambda t}
			\quad\quad
			\forall t\geq T_{*}.
			\label{th:ODE-slow-lambda}
		\end{equation}
		
	\end{enumerate}
\end{prop}

As a consequence, if we want \emph{all} solutions to (\ref{eqn:ODE})
to decay faster than a given exponential, we are forced to choose a
damping coefficient $\delta(t)$ which alternates intervals where it is
large enough, and intervals where it is smaller than $\lambda$.  If we
want solutions to decay faster than all exponentials, we also need
$\delta(t)$ to be unbounded.

In some sense these are the unique essential features.  In the
following result we show that any fixed exponential decay rate can be
achieved through a periodic damping coefficient which is greater than
or equal to $\lambda-\ep$, and has Lipschitz constant equal to $\ep$,
and hence it exhibits very slow transitions from small to large values.

\begin{thm}[Single ODE, with borderline constraints on $\delta(t)$]\label{thm:ODE-fixed-lip}
	Let $\lambda$ and $R$ be positive real numbers, and let 
	$\ep\in(0,\lambda)$. Then there exist a positive real number 
	$t_{0}$, with
	\begin{equation}
		t_{0}\leq 
		16\left(\frac{\pi}{\ep}+1\right)R+
		\frac{2(\pi+1)}{\ep}+\frac{8}{\lambda}+
		2+8\log 2,
		\label{th:est-t0-lip}
	\end{equation}
	and a $t_{0}$-periodic function $\delta:[0,+\infty)\to[0,+\infty)$
	such that
	$$|\delta(t)-\delta(s)|\leq\ep|t-s|
	\quad\quad
	\forall t\geq 0\quad\forall s\geq 0,$$
	$$\delta(t)\geq\lambda-\ep
	\quad\quad
	\forall t\geq 0,$$
	and such that every solution to (\ref{eqn:ODE}) satisfies
	\begin{equation}
		|u'(t)|^{2}+\lambda^{2}|u(t)|^{2}\leq
		\left(|u'(0)|^{2}+\lambda^{2}|u(0)|^{2}\right)
		\exp\left(-R(t-t_{0})^{+}\right)
		\quad\quad
		\forall t\geq 0.
		\label{th:ODE-fixed-lip}
	\end{equation}
\end{thm}

Theorem~\ref{thm:ODE-fixed-lip} deals with the case of a single 
ordinary differential equation, but an analogous result holds true 
also for systems or the abstract equation (\ref{eqn:PDE}). In order to 
contain this paper in a reasonable length we spare the reader from 
the details.

More delicate is achieving decay rates faster than all exponentials
through damping coefficients with small Lipschitz constant.  This
could be done at most in the same spirit of Theorem~\ref{thm:PDE-top},
the reason being as usual that now $t_{0}$ depends on $R$.  We do not
address this issue in this paper.

\setcounter{equation}{0}
\section{Heuristics}\label{sec:heuristics}

In this section we provide an informal description of the strategy of 
our proofs. Our aim is clarifying why pulsating damping coefficients 
are more effective when we are interested in damping all solutions to 
an equation or system.

Let us start with the single ordinary differential equation
(\ref{eqn:ODE}).  We want to design $\delta(t)$ so that all solutions
decay as fast as possible.  The first naive idea is to choose
$\delta(t)$ very large.  Bending the rules a little bit, we can even
imagine to choose a damping coefficient which is not a function, but
a Dirac delta function (which is actually a measure) concentrated at
time $t=0$, or even better a delta function multiplied by a large
enough constant $k$.

An easy calculation shows that such an extreme damping has a great
effect on the solution with initial data $u(0)=0$ and $u'(0)=1$, whose
energy is instantly reduced by a factor $e^{-k}$.  On the contrary, it
has no effect on the orthogonal solution with initial data $u(0)=1$
and $u'(0)=0$.  This can be explained by observing that the damping
coefficient multiplies $u'(t)$, and in the case of the second solution
this time-derivative vanishes when the delta function acts.  The
effect on any other solution is a linear combination of the two,
namely highly reducing the time-derivative but leaving the function
untouched.

This apparently inconclusive approach suggests a first strategy: if we
want to dampen a single solution, we can use a delta function acting when
the energy of the solution is concentrated on the time-derivative.  So
we take again the solution with initial data $u(0)=1$ and $u'(0)=0$,
and we apply no dissipation until $u(t)$ vanish.  This happens for the
first time when $t=\pi/(2\lambda)$.  At that point we apply a second
delta function.

Summing up, a first delta function at time $t=0$ cuts the first
solution, then the damping coefficient vanishes until the second delta
function at time $t=\pi/(2\lambda)$ cuts the second solution.  What
happens to all other solutions?  Since the equation is linear, cutting
two linearly independent solutions is equivalent to cutting all
solutions.  If we repeat this procedure by periodicity, we can achieve
any exponential decay rate.  If at each reiteration we increase the
multiplicative constant in front of the delta functions, we can
achieve any given decay rate.  This is the idea behind the proofs of
the results stated in section~\ref{sec:ODE}, and the point where the
special time (\ref{defn:t0-ODE}) and the profiles of
Figure~\ref{fig:delta-ODE} come into play.

The idea for systems is a simple generalization. We use a first block 
of two delta functions with time-gap of $\pi/(2\lambda_{1})$ in order 
to dampen solutions of the first equation, then we use two more delta 
functions with time-gap of $\pi/(2\lambda_{2})$ in order to dampen 
solutions of the second equation, and so on. In other words, we take 
care of the equations of the system one by one. We end up with the 
special time (\ref{defn:t0-system}) and the profiles of 
Figure~\ref{fig:delta-system}. We are quite skeptic about the 
possibility of reducing the time (\ref{defn:t0-system}), unless 
$\lambda_{i}$'s satisfy special rationality conditions.

When we deal with a partial differential equation, which we regard as
a system of countably many ordinary differential equations, we exploit
a mixed strategy.  If we want to achieve a given exponential decay
rate, a suitable constant damping coefficient does the job for all
components corresponding to large enough eigenvalues.  Thus we are
left with cutting a finite number of components, and this can be done
as in the case of finite systems.  As a consequence, now a good
damping coefficient consists in a constant damping half the time,
alternated with a train of delta functions in the remaining half of
the time.  This is the idea behind the proof of
Theorem~\ref{thm:PDE-fixed} and the profile of
Figure~\ref{fig:delta-PDE}.

In the case of partial differential equations, things are more complex
if we want to achieve a decay rate faster than all exponentials.
Indeed, when we reiterate the procedure, a better decay requires more
components to be treated separately, and in turn this implies a longer
wait.  The compromise between faster decay rates and longer waiting
times gives rise to the operator dependent rate $\varphi(t)$ of
Theorem~\ref{thm:PDE-top}.

This discussion motivates also the last part of the paper.  Indeed, a
train of delta functions (or a suitable approximation) emerged as a
common pattern of the damping coefficients which realize fast decay
rates.  In a first stage this led us to suspect that a bound on the
Lipschitz constant of the coefficient, or the impossibility to attain
values close to zero, could yield a bound from below on the decay rate
of solutions.

In Theorem~\ref{thm:ODE-fixed-lip} we show that this is not the case,
because any exponential decay rate can be realized through a damping
coefficient $\delta(t)$ with arbitrarily small Lipschitz constant.  The
construction of $\delta(t)$ is more involved, but once
again it acts in two steps.  In a first phase the coefficient grows
and kills a first solution.  Then the coefficients goes below $\lambda$
and stays there until the orthogonal solution has rotated enough so
that it is ready to be damped by a second growth of the coefficient.

This two-phase action (destroy the first solution, wait for rotation, 
destroy the second solution) seems to be the quintessence of all the 
story.

\setcounter{equation}{0}
\section{Proofs}\label{sec:proofs}

In this section we prove our main results, following the same scheme
of the statement section.  We begin by investigating how solutions to
(\ref{eqn:PDE}) depend on the damping coefficient.

\begin{lemma}[Continuous dependence on the damping coefficient]\label{lemma:approx}
	Let $H$ be a Hilbert space, let $A$ be a self-adjoint nonnegative
	linear operator on $H$ with dense domain $D(A)$, and let $T$ be a
	positive real number.  Let $\delta_{1}:[0,T]\to[0,+\infty)$ and
	$\delta_{2}:[0,T]\to[0,+\infty)$ be two bounded measurable
	functions.  Let $u_{1}(t)$ and $u_{2}(t)$ be the solutions to
	(\ref{eqn:PDE}) with $\delta(t)$ replaced by $\delta_{1}(t)$ and
	$\delta_{2}(t)$, respectively, and with initial data
	$$u_{1}(0)=u_{2}(0)=u_{0}\in D(A^{1/2}), \hspace{3em}
	u_{1}'(0)=u_{2}'(0)=u_{1}\in H.$$
	
	Then for every $t\in[0,T]$ the following estimate holds true
	\begin{eqnarray}
		|u_{2}'(t)-u_{1}'(t)|^{2}+|A^{1/2}(u_{2}(t)-u_{1}(t))|^{2} & \leq & 
		2\left(|u_{1}|^{2}+|A^{1/2}u_{0}|^{2}\right)e^{2t}\cdot
		\nonumber  \\
		 &  & \mbox{}\cdot \int_{0}^{t}
		|\delta_{2}(s)-\delta_{1}(s)|^{2}\,ds.
		\label{th:u2-u1}
	\end{eqnarray}
\end{lemma}

\paragraph{\textmd{\textit{Proof}}}

To begin with, we observe that
\begin{equation}
	|u_{1}'(t)|^{2}+|A^{1/2}u_{1}(t)|^{2}\leq
	|u_{1}|^{2}+|A^{1/2}u_{0}|^{2}
	\quad\quad
	\forall t\in[0,T].
	\label{energy-u1}
\end{equation}

This inequality holds true because the left-hand side is a
nonincreasing function of time.  Now let us set
$$E(t):=|u_{2}'(t)-u_{1}'(t)|^{2}+
|A^{1/2}(u_{2}(t)-u_{1}(t))|^{2}.$$

An easy computation shows that
$$E'(t)=-4\delta_{2}(t)|u_{2}'(t)-u_{1}'(t)|^{2}+
4(\delta_{1}(t))-\delta_{2}(t))\cdot 
\langle u_{1}'(t),u_{2}'(t)-u_{1}'(t)\rangle.$$

The first term in the right-hand side is less than or equal to zero. 
Keeping (\ref{energy-u1}) into account, we can estimate the second 
term and obtain that
\begin{eqnarray*}
	E'(t) & \leq & 2|u_{2}'(t)-u_{1}'(t)|^{2}+
	2|u_{1}'(t)|^{2}\cdot|\delta_{2}(t)-\delta_{1}(t)|^{2} \\
	\noalign{\vspace{1ex}}
	 & \leq & 2E(t)+
	 2\left(|u_{1}|^{2}+|A^{1/2}u_{0}|^{2}\right)
	 |\delta_{2}(t)-\delta_{1}(t)|^{2}.
\end{eqnarray*}

Integrating this differential inequality, and recalling that $E(0)=0$ 
because the initial conditions of $u_{2}(t)$ and $u_{1}(t)$ are the 
same, we conclude that
$$E(t)\leq 2\left(|u_{1}|^{2}+|A^{1/2}u_{0}|^{2}\right)
e^{2t}\int_{0}^{t}
|\delta_{2}(s)-\delta_{1}(s)|^{2}\,ds
\quad\quad
\forall t\in[0,T],$$
which proves (\ref{th:u2-u1}).\qed

\subsection{Ordinary differential equations}\label{sec:proofs-ODE}

The following result is the fundamental tool in our theory.

\begin{lemma}[Decay for two orthogonal solutions to a single ODE]\label{lemma:basic-vw}
	Let $\lambda$ and $M$ be positive real numbers, and let $t_{0}$ 
	be defined by (\ref{defn:t0-ODE}).  For every positive integer 
	$n$, let us consider the function 
	$\delta_{n}:[0,t_{0}]\to[0,+\infty)$ defined by
	\begin{equation}
		\delta_{n}(t):=\left\{
		\begin{array}{ll}
			Mn\quad & \mbox{if }t\in [0,1/n]\cup[t_{0}-1/n,t_{0}],  \\
			\noalign{\vspace{0.5ex}}
			0 & \mbox{otherwise},
		\end{array}
		\right.
		\label{defn:delta-n}
	\end{equation}
	and the differential equation
	\begin{equation}
		u''(t)+2\delta_{n}(t)u'(t)+\lambda^{2}u(t)=0.
		\label{ODE:delta-n}
	\end{equation}
	
	Let $v_{n}(t)$ be the solution with initial data $v_{n}(0)=0$ and
	$v_{n}'(0)=1$.  Let $w_{n}(t)$ be the solution with initial data
	$w_{n}(0)=1/\lambda$ and $w_{n}'(0)=0$.
	
	Then
	\begin{equation}
		\limsup_{n\to +\infty}
		\left(|v_{n}'(t_{0})|^{2}+\lambda^{2}|v_{n}(t_{0})|^{2}\right)
		\leq e^{-4M},
		\label{th:lim-vn}
	\end{equation}
	\begin{equation}
		\limsup_{n\to +\infty}
		(|w_{n}'(t_{0})|^{2}+\lambda^{2}|w_{n}(t_{0})|^{2})
		\leq e^{-4M}.
		\label{th:lim-wn}
	\end{equation}
\end{lemma}

\paragraph{\textmd{\textit{Proof}}}

For every solution $u(t)$ to (\ref{ODE:delta-n}), let us consider its 
energy
$$E_{u}(t):=|u'(t)|^{2}+\lambda^{2}|u(t)|^{2}.$$

A simple computation of the time-derivative shows that $E_{u}(t)$ is 
nonincreasing. Let us set for simplicity
$$t_{n}:=1/n,
\hspace{4em}
s_{n}:=t_{0}-1/n,$$
and let us assume that $n$ is large enough to that $t_{n}<s_{n}$, 
and hence the two intervals where $\delta_{n}(t)=Mn$ are disjoint.

\subparagraph{\textmd{\textit{Estimate on $v_{n}(t)$}}}

Since $E_{v_{n}}(t)\leq E_{v_{n}}(0)=1$, it follows that 
$|v_{n}'(t)|\leq 1$ for every $t\in[0,t_{0}]$. Thus from the mean 
value theorem we obtain that $|v_{n}(0)-v_{n}(t_{n})|\leq t_{n}$,
and hence
\begin{equation}
	\lim_{n\to +\infty}v_{n}(t_{n})=0.
	\label{lim-vn-tn}
\end{equation}

Let us consider now the time-derivative. To this end, we interpret 
(\ref{ODE:delta-n}) as a first order linear equation in $u'(t)$, with 
forcing term $-\lambda^{2}u(t)$. Integrating this differential 
equation we obtain that
$$v_{n}'(t)=v_{n}'(0)e^{-2Mnt}+\int_{0}^{t}
e^{2Mn(s-t)}\lambda^{2}v_{n}(s)\,ds
\quad\quad
\forall t\in[0,t_{n}].$$

Now we set $t=t_{n}$, we recall that $v_{n}'(0)=1$, and we pass to the
limit as $n\to +\infty$.  The integrand is bounded because $s\leq t$
and $|v_{n}(s)|$ is bounded owing to the energy estimate. Since 
$t_{n}\to 0$  the integral tends to 0, and hence
\begin{equation}
	\lim_{n\to +\infty}v_{n}'(t_{n})=e^{-2M}.
	\label{lim-vn'-tn}
\end{equation}

From (\ref{lim-vn-tn}) and (\ref{lim-vn'-tn}) we conclude that
$$\lim_{n\to +\infty}
\left(|v_{n}'(t_{n})|^{2}+\lambda^{2}|v_{n}(t_{n})|^{2}\right)
=e^{-4M},$$
which implies (\ref{th:lim-vn}) because the energy is nonincreasing 
with time.

\subparagraph{\textmd{\textit{Estimate on $w_{n}(t)$}}}

Let us begin with the interval $[0,t_{n}]$. The same argument 
exploited in the case of $v_{n}(t)$ now leads to
\begin{equation}
	\lim_{n\to +\infty}w_{n}(t_{n})=1/\lambda,
	\hspace{4em}
	\lim_{n\to +\infty}w_{n}'(t_{n})=0.
	\label{lim-wn-tn}
\end{equation}

Roughly speaking, this means that nothing changes for $w_{n}(t)$ in 
the interval $[0,t_{n}]$.

Let us consider now the interval $(t_{n},s_{n})$, where 
$\delta_{n}(t)$ is identically 0. An easy computation shows that in 
this interval $w_{n}(t)$ is given by the explicit formula
$$w_{n}(t)=w_{n}(t_{n})\cos(\lambda(t-t_{n}))+
\frac{w_{n}'(t_{n})}{\lambda}\sin(\lambda(t-t_{n})).$$

Setting $t=s_{n}$ we find that
$$w_{n}(s_{n})=
w_{n}(t_{n})\cos\left(\frac{\pi}{2}-\frac{2\lambda}{n}\right)+
\frac{w_{n}'(t_{n})}{\lambda}
\sin\left(\frac{\pi}{2}-\frac{2\lambda}{n}\right),$$
$$w_{n}'(s_{n})=-\lambda w_{n}(t_{n})
\sin\left(\frac{\pi}{2}-\frac{2\lambda}{n}\right)+
w_{n}'(t_{n})\cos\left(\frac{\pi}{2}-\frac{2\lambda}{n}\right).$$

Passing to the limit as $n\to +\infty$, and keeping (\ref{lim-wn-tn}) 
into account, we deduce that
$$\lim_{n\to +\infty}w_{n}(s_{n})=0,
\hspace{4em}
\lim_{n\to +\infty}w_{n}'(s_{n})=-1.$$

Roughly speaking, this means that the interval $[t_{n},s_{n}]$ has
produced a rotation of $w_{n}(t)$ in the phase space, with the effect
of moving all the energy on the derivative.

Let us finally consider the interval $[s_{n},t_{0}]$, where we argue 
as we did in $[0,t_{n}]$ with the function $v_{n}(t)$. Due to the  
uniform bound on $w_{n}'(t)$ coming from the energy estimate, from 
the mean value theorem we deduce that
\begin{equation}
	\lim_{n\to +\infty}w_{n}(t_{0})=
	\lim_{n\to +\infty}w_{n}(s_{n})=0.
	\label{lim-wn-T0}
\end{equation}

As for the derivative, once again we interpret (\ref{ODE:delta-n}) as 
a first order linear equation in $u'(t)$, and we find that
$$w_{n}'(t)=w_{n}'(s_{n})e^{-2Mn(t-s_{n})}+\int_{s_{n}}^{t}
e^{2Mn(s-t)}\lambda^{2}w_{n}(s)\,ds
\quad\quad
\forall t\in[s_{n},t_{0}].$$

Once again the integrand is bounded because $s\leq t$ and $w_{n}(s)$ 
is uniformly bounded owing to the energy estimate. Setting $t=t_{0}$, 
and passing to the limit as $n\to +\infty$, we conclude that
\begin{equation}
	\lim_{n\to +\infty}w_{n}'(t_{0})=-e^{-2M}.
	\label{lim-wn'-T0}
\end{equation}

From (\ref{lim-wn-T0}) and (\ref{lim-wn'-T0}) we deduce 
(\ref{th:lim-wn}).\qed
\medskip

A careful inspection of the proof reveals that the inequality in
(\ref{th:lim-wn}) is actually an equality, and the limsup is actually
a limit. With similar arguments one could show that the same is true 
in (\ref{th:lim-vn}) (it is enough to follow the solution $v_{n}(t)$ 
until the end of the interval). In any case, (\ref{th:lim-vn}) and 
(\ref{th:lim-wn}) are what we need in the sequel.

In the next result we apply Lemma~\ref{lemma:basic-vw} in order to 
dampen all solutions to~(\ref{eqn:ODE}).

\begin{lemma}[Decay for all solutions to a single ODE]\label{lemma:basic-u}
	Let $\lambda$ and $M$ be positive real numbers, and let $t_{0}$ 
	be defined by (\ref{defn:t0-ODE}). 
	
	Then there exists a bounded measurable damping coefficient
	$\delta:[0,t_{0}]\to[0,+\infty)$ such that every solution to
	(\ref{eqn:ODE}) satisfies
	\begin{equation}
		|u'(t_{0})|^{2}+\lambda^{2}|u(t_{0})|^{2}\leq
		\left(|u'(0)|^{2}+\lambda^{2}|u(0)|^{2}\right)
		\cdot 2e^{-M}.
		\label{th:u-T0}
	\end{equation}
	
	Furthermore, one can choose $\delta(t)$ such that
	\begin{itemize}
		\item  either it is of the form (\ref{defn:delta-n}) for a 
		large enough $n$,
	
		\item  or it is of class $C^{\infty}$ and all its 
		time-derivatives vanish both in $t=0$ and in $t=t_{0}$.
	\end{itemize}
\end{lemma}

\paragraph{\textmd{\textit{Proof}}}

We claim that (\ref{th:u-T0}) holds true for all solutions to
(\ref{eqn:ODE}) if the damping coefficient is of the form
(\ref{defn:delta-n}) with $n$ large enough.  To this end, let
$\delta_{n}(t)$, $v_{n}(t)$ and $w_{n}(t)$ be defined as in
Lemma~\ref{lemma:basic-vw}.  Due to (\ref{th:lim-vn}) and
(\ref{th:lim-wn}), the two inequalities
\begin{equation}
	|v_{n}'(t_{0})|^{2}+\lambda^{2}|v_{n}(t_{0})|^{2}
	\leq e^{-M},
	\hspace{3em}
	|w_{n}'(t_{0})|^{2}+\lambda^{2}|w_{n}(t_{0})|^{2}
	\leq e^{-M}
	\label{vn-wn-eM}
\end{equation}
hold true provided that $n$ is large enough. Every solution $u(t)$ to 
equation (\ref{eqn:ODE}) with $\delta(t):=\delta_{n}(t)$
is a linear combination of $v_{n}(t)$ and $w_{n}(t)$, and more
precisely
$$u(t)=u'(0)v_{n}(t)+\lambda u(0)w_{n}(t).$$

It follows that
\begin{eqnarray*}
	|u'(t_{0})|^{2}+\lambda^{2}|u(t_{0})|^{2} & \leq & 
	2|u'(0)|^{2}\cdot|v_{n}'(t_{0})|^{2}+
	2\lambda^{2}|u(0)|^{2}\cdot|w_{n}'(t_{0})|^{2}   \\
	\noalign{\vspace{1ex}}
	 &  & \mbox{}+2\lambda^{2}|u'(0)|^{2}\cdot|v_{n}(t_{0})|^{2}+
	2\lambda^{4}|u(0)|^{2}\cdot|w_{n}(t_{0})|^{2}  \\
	\noalign{\vspace{1ex}}
	 & = & 2|u'(0)|^{2}\cdot\left(
	 |v_{n}'(t_{0})|^{2}+\lambda^{2}|v_{n}(t_{0})|^{2}
	 \right)   \\
	 \noalign{\vspace{1ex}}
	 & & \mbox{}+2\lambda^{2}|u(0)|^{2}\cdot\left(
	 |w_{n}'(t_{0})|^{2}+\lambda^{2}|w_{n}(t_{0})|^{2}
	 \right),
\end{eqnarray*}
so that (\ref{th:u-T0}) follows from (\ref{vn-wn-eM}).  This proves 
our original claim.

If we want a smooth damping coefficient, we need a two steps 
approximation. First of all we choose $n_{0}\in\n$ such that
$$|v_{n_{0}}'(t_{0})|^{2}+\lambda^{2}|v_{n_{0}}(t_{0})|^{2}\leq 
e^{-2M},
\hspace{3em}
|w_{n_{0}}'(t_{0})|^{2}+\lambda^{2}|w_{n_{0}}(t_{0})|^{2}
\leq e^{-2M}.$$

Then we approximate $\delta_{n_{0}}(t)$ with a nonnegative function
$\delta(t)$ of class $C^{\infty}$ with the property that all its
time-derivatives vanish both in $t=0$ and in $t=t_{0}$.  Let $v(t)$
and $w(t)$ denote the corresponding solutions to (\ref{eqn:ODE}) with
the same initial data of $v_{n_{0}}(t)$ and $w_{n_{0}}(t)$,
respectively.  From Lemma~\ref{lemma:approx} we deduce that, if
$\delta(t)$ is close enough to $\delta_{n_{0}}(t)$ in
$L^{2}((0,t_{0}))$, then
$$|v'(t_{0})|^{2}+\lambda^{2}|v(t_{0})|^{2}\leq e^{-M}, \hspace{3em}
|w'(t_{0})|^{2}+\lambda^{2}|w(t_{0})|^{2}\leq e^{-M}.$$

Since every solution to (\ref{eqn:ODE}) is a linear combination of 
$v(t)$ and $w(t)$, we can conclude exactly as before.\qed

\subsubsection*{Proof of Theorem~\ref{thm:ODE-fixed}}

Let $M$ be such that $2e^{-M}=e^{-Rt_{0}}$.  Let us consider any
function $\delta(t)$ provided by Lemma~\ref{lemma:basic-u} with this
choice of $M$, and let us extend it by periodicity to the half line
$t\geq 0$.  Let $u(t)$ be any corresponding solution to
(\ref{eqn:ODE}), and let us consider the usual energy
$E(t):=|u'(t)|^{2}+\lambda^{2}|u(t)|^{2}$.

The estimate of Lemma~\ref{lemma:basic-u} can be applied in all 
intervals of the form $[kt_{0},(k+1)t_{0}]$, yielding 
that
$$E((k+1)t_{0})\leq E(kt_{0})\cdot 2e^{-M}=
E(kt_{0})\cdot e^{-Rt_{0}}
\quad\quad
\forall k\in\n.$$

Therefore, an easy induction proves that
$$E(kt_{0})\leq E(0)\cdot e^{-kRt_{0}}
\quad\quad
\forall k\in\n.$$

Since $E(t)$ is nonincreasing, this implies that
$$E(t)\leq E(0)\exp\left(-R(t-t_{0})^{+}\right)
\quad\quad
\forall t\geq 0,$$
which is equivalent to (\ref{th:ODE-fixed}).

If we want $\delta(t)$ to be piecewise constant, it is enough to take a
function with this property from Lemma~\ref{lemma:basic-u}.  If we
want $\delta(t)$ to be of class $C^{\infty}$, it is enough to take
from Lemma~\ref{lemma:basic-u} a function of class $C^{\infty}$ whose
time-derivatives of any order vanish at the endpoints of the
interval.  This condition guarantees that the periodic extension
remains of class $C^{\infty}$.\qed

\subsubsection*{Proof of Theorem~\ref{thm:ODE-any}}

Let $M_{k}$ be a sequence of positive real numbers such that
\begin{equation}
	2e^{-M_{0}}\leq\varphi(2t_{0}),
	\label{defn:M0}
\end{equation}
and
\begin{equation}
	\varphi((k+1)t_{0})\cdot 2e^{-M_{k}}\leq
	\varphi((k+2)t_{0})
	\quad\quad
	\forall k\geq 1.
	\label{defn:Mk}
\end{equation}

For every $k\in\n$, let $\delta_{k}:[0,t_{0}]\to[0,+\infty)$ be one 
of the functions provided by Lemma~\ref{lemma:basic-u}, applied with 
$M:=M_{k}$. 
Let us define $\delta:[0,+\infty)\to[0,+\infty)$ by 
glueing together all these functions, namely by setting
$$\delta(t):=\delta_{k}(t-kt_{0})
\quad\quad
\forall t\in[kt_{0},(k+1)t_{0}).$$

Let $u(t)$ be any corresponding solution to (\ref{eqn:ODE}), and let
us consider the usual energy
$$E(t):=|u'(t)|^{2}+\lambda^{2}|u(t)|^{2}.$$

We claim that $E(t)\leq E(0)\varphi(t)$ for all $t\geq t_{0}$, which 
is equivalent to (\ref{th:ODE-any}). In order to prove this result, 
it is enough to show that
\begin{equation}
	E(kt_{0})\leq E(0)\cdot\varphi((k+1)t_{0})
	\quad\quad
	\forall k\geq 1.
	\label{th:ODE-any-k}
\end{equation}

Indeed, since both $E(t)$ and $\varphi(t)$ are nonincreasing, when 
$t\in[kt_{0},(k+1)t_{0}]$ for some $k\geq 1$ it turns out that
$$E(t)\leq
E(kt_{0})\leq
E(0)\cdot\varphi((k+1)t_{0})\leq
E(0)\cdot\varphi(t).$$

In order to prove (\ref{th:ODE-any-k}), we repeatedly apply 
Lemma~\ref{lemma:basic-u}. Since in the interval $[0,t_{0}]$ the function 
$\delta(t)$ coincides with $\delta_{0}(t)$, from 
Lemma~\ref{lemma:basic-u} and (\ref{defn:M0}) we deduce that
$$E(t_{0})\leq E(0)\cdot 2e^{-M_{0}}\leq E(0)\cdot\varphi(2t_{0}),$$
which proves (\ref{th:ODE-any-k}) in the case $k=1$.

Now we proceed by induction. Let us assume that (\ref{th:ODE-any-k}) 
holds true for some positive integer $k$. In the interval 
$[kt_{0},(k+1)t_{0}]$ the function $\delta(t)$ is a time 
translation of $\delta_{k}(t)$, and thus we can apply 
Lemma~\ref{lemma:basic-u} up to this time translation. Keeping 
(\ref{defn:Mk}) into account, we deduce that
$$E((k+1)t_{0})\leq E(kt_{0})\cdot 2e^{-M_{k}}\leq
E(0)\cdot\varphi((k+1)t_{0})\cdot 2e^{-M_{k}}\leq
E(0)\cdot\varphi((k+2)t_{0}),$$
which completes the induction.

The possibility of choosing a damping coefficient satisfying further 
requirements depends on the analogous possibility in 
Lemma~\ref{lemma:basic-u}\qed

\subsection{Systems of ordinary differential equations}\label{sec:proofs-system}

The following result is the natural generalization of 
Lemma~\ref{lemma:basic-u} to systems.

\begin{lemma}[Decay for all solutions to a system]\label{lemma:basic-u-system}
	Let $n$ be a positive integer, let
	$(\lambda_{1},\ldots,\lambda_{n})\in(0,+\infty)^{n}$, let $M$ be a
	positive real number, and let $t_{0}$ be defined by
	(\ref{defn:t0-system}).
	
	Then there exists a bounded measurable damping coefficient
	$\delta:[0,t_{0}]\to[0,+\infty)$ such that every solution to
	(\ref{eqn:ODE-system}) satisfies
	\begin{equation}
		\sum_{k=1}^{n}
		\left(|u_{k}'(t_{0})|^{2}+
		\lambda_{k}^{2}|u_{k}(t_{0})|^{2}\right)\leq
		\left[\sum_{k=1}^{n}
		\left(|u_{k}'(0)|^{2}+
		\lambda_{k}^{2}|u_{k}(0)|^{2}\right)\right]
		\cdot 2e^{-M}.
		\label{th:u-T0-k}
	\end{equation}
	
	Furthermore, one can choose $\delta(t)$ such that
	\begin{itemize}
		\item  either it has the profile shown in 
		Figure~\ref{fig:delta-system} on the left,
	
		\item  or it is of class $C^{\infty}$ and all its 
		time-derivatives vanish both in $t=0$ and in $t=t_{0}$.
	\end{itemize}
\end{lemma}

\paragraph{\textmd{\textit{Proof}}}

Let us apply Lemma~\ref{lemma:basic-u} to the $k$-th equation of the 
system. Setting $t_{0k}:=\pi/(2\lambda_{k})$, we obtain a function 
$\delta_{k}:[0,t_{0k}]\to[0,+\infty)$ such that every solution to the 
$k$-th equation of the system satisfies
$$|u_{k}'(t_{0k})|^{2}+\lambda_{k}^{2}|u_{k}(t_{0k})|^{2}\leq
\left(|u_{k}'(0)|^{2}+\lambda_{k}^{2}|u_{k}(0)|^{2}\right)
\cdot 2e^{-M}.$$

Now we observe that $t_{0}=t_{01}+\ldots+t_{0n}$, and we define 
$\delta:[0,t_{0}]\to[0,+\infty)$ by glueing together the functions 
$\delta_{k}(t)$ defined above. More precisely, we partition 
$[0,t_{0}]$ into $n$ subintervals $[s_{k-1},s_{k}]$ with $s_{0}=0$, 
$s_{n}=t_{0}$ and $s_{k}-s_{k-1}=t_{0k}$, and then we set
$$\delta(t):=\delta_{k}(t-s_{k-1})
\quad\quad
\forall t\in[s_{k-1},s_{k}).$$

Let $(u_{1}(t),\ldots,u_{n}(t))$ be a corresponding solution to the 
system (\ref{eqn:ODE-system}), and let us set
$$E_{k}(t):=|u_{k}'(t)|^{2}+\lambda_{k}^{2}|u_{k}(t)|^{2}.$$

We claim that the energy $E_{k}(t)$ of the $k$-th component is 
reduced by a factor $2e^{-M}$ in the interval $[s_{k-1},s_{k}]$. 
Indeed in this interval $\delta(t)$ coincides with $\delta_{k}(t)$ up 
to a time translation, hence from Lemma~\ref{lemma:basic-u}
we deduce that
$$E_{k}(t_{0})\leq
E_{k}(s_{k})\leq
E_{k}(s_{k-1})\cdot 2e^{-M}\leq
E_{k}(0)\cdot 2e^{-M}.$$

Summing over all indices $k$ from 1 to $n$ we obtain
(\ref{th:u-T0-k}).

If we want a damping coefficient with the profile shown in
Figure~\ref{fig:delta-system} on the left, it is enough that all 
functions $\delta_{k}(t)$ are of the form (\ref{defn:delta-n}), with 
the same value of $n$ for all $k$'s (this is possible 
provided that $n$ is large enough).

If we want a damping coefficient of class $C^{\infty}$, it is enough to
choose all functions $\delta_{k}(t)$ of class $C^{\infty}$ with
all time-derivatives vanishing at the endpoints of the interval.  This
condition guarantees that the glueing procedure yields a function
which is still of class $C^{\infty}$.\qed

\subsubsection*{Proof of Theorem~\ref{thm:system-fixed} and
Theorem~\ref{thm:system-any}}

The arguments are analogous to the proofs of
Theorem~\ref{thm:ODE-fixed} and Theorem~\ref{thm:ODE-any}, just 
starting from Lemma~\ref{lemma:basic-u-system} instead of 
Lemma~\ref{lemma:basic-u}.\qed

\begin{rmk}\label{rmk:eigenspaces}
	\begin{em}
		As already mentioned in the introduction, now it should be
		clear from the proofs that the conclusions of
		Theorem~\ref{thm:system-fixed} and
		Theorem~\ref{thm:system-any} hold true for the general
		equation (\ref{eqn:PDE}), provided that the spectrum of $A$ is
		finite, even if the dimension of eigenspaces is infinite.
	\end{em}
\end{rmk}

\subsection{Partial differential equations}\label{sec:proofs-PDE}

When the operator $A$ has infinitely many eigenvalues, a constant 
dissipation is enough to dampen all components corresponding to large 
enough eigenvalues. This is the content of next result.

\begin{lemma}[PDE with constant dissipation]\label{lemma:PDE-coercive}
	Let $H$ be a Hilbert space, and let $A$ be a self-adjoint
	nonnegative operator on $H$ with dense domain $D(A)$. Let $M$ be 
	a positive real number, and let us assume that
	\begin{equation}
		|A^{1/2}x|^{2}\geq 2M^{2}|x|^{2}
		\quad\quad
		\forall x\in D(A^{1/2}).
		\label{hp:A-coercive}
	\end{equation}
	
	Then every weak solution $u(t)$ to
	$$u''(t)+2Mu'(t)+Au(t)=0$$
	satisfies
	\begin{equation}
		|u'(t)|^{2}+|A^{1/2}u(t)|^{2}\leq
		\left(|u'(0)|^{2}+|A^{1/2}u(0)|^{2}\right)
		\cdot 8e^{-2Mt}
		\quad\quad
		\forall t\geq 0.
		\label{th:PDE-coercive}
	\end{equation}
\end{lemma}

\paragraph{\textmd{\textit{Proof}}}

Let us consider the usual energy
$$E(t):=|u'(t)|^{2}+|A^{1/2}u(t)|^{2},$$
and the modified energy
$$\Ehat(t):=|u'(t)|^{2}+|A^{1/2}u(t)|^{2}+
2M\langle u'(t),u(t)\rangle.$$

An easy calculation shows that
\begin{equation}
	\Ehat'(t)=-2M\Ehat(t)
	\quad\quad
	\forall t\geq 0.
	\label{eqn:deriv-Ehat}
\end{equation}

We claim that
\begin{equation}
	\frac{1}{4}E(t)\leq\Ehat(t)\leq 2E(t)
	\quad\quad
	\forall t\geq 0.
	\label{E-Ehat}
\end{equation}

Indeed, from the inequality
$$2M\langle u'(t),u(t)\rangle\leq
\frac{1}{2}|u'(t)|^{2}+2M^{2}|u(t)|^{2}$$
and assumption
(\ref{hp:A-coercive}) it follows that
$$\Ehat(t)\leq E(t)+\frac{1}{2}|u'(t)|^{2}+|A^{1/2}u(t)|^{2}
\leq 2E(t).$$

On the other hand, from the inequality
$$2M\langle u'(t),u(t)\rangle\geq
-\frac{3}{4}|u'(t)|^{2}-\frac{4}{3}M^{2}|u(t)|^{2}$$
and assumption
(\ref{hp:A-coercive}) it follows that
$$\Ehat(t)\geq E(t)-\frac{3}{4}|u'(t)|^{2}-\frac{2}{3}|A^{1/2}u(t)|^{2}
\geq\frac{1}{4}E(t).$$

From (\ref{eqn:deriv-Ehat}) and (\ref{E-Ehat}) we conclude that
$$E(t)\leq 4\Ehat(t)=4\Ehat(0)e^{-2Mt}\leq 8E(0)e^{-2Mt},$$
which is equivalent to (\ref{th:PDE-coercive}).
\qed

The following result is the generalization of 
Lemma~\ref{lemma:basic-u} and Lemma~\ref{lemma:basic-u-system} to the 
infinite dimensional setting.

\begin{lemma}[Decay for all solutions to a PDE]\label{lemma:basic-PDE}
	Let $H$, $A$, $\lambda_{k}$ be as in Theorem~\ref{thm:PDE-fixed}. 
	Let $R$ be a positive real number, and let 
	\begin{equation}
		T_{R}:=\frac{\pi}{2}
		\sum_{\lambda_{k}^{2}\leq 2(R+\lambda_{1})^{2}}
		\frac{1}{\lambda_{k}}.
		\label{defn:TR}
	\end{equation}
	
	Then there exists a bounded measurable damping coefficient
	$\delta:[0,2T_{R}]\to[0,+\infty)$ such that every weak solution to
	(\ref{eqn:PDE}) satisfies
	\begin{equation}
		|u'(2T_{R})|^{2}+|A^{1/2}u(2T_{R})|^{2}\leq
		\left(|u'(0)|^{2}+|A^{1/2}u(0)|^{2}\right)\cdot
		e^{-R\cdot 2T_{R}}.
		\label{th:u-2TR}
	\end{equation}
	
	Furthermore, one can choose $\delta(t)$ such that
	\begin{itemize}
		\item either it has the profile shown in 
		Figure~\ref{fig:delta-PDE},
	
		\item  or it is of class $C^{\infty}$ and all its 
		time-derivatives vanish both in $t=0$ and in $t=2T_{R}$.
	\end{itemize}
\end{lemma}

\paragraph{\textmd{\textit{Proof}}}

Let us write $H$ as a direct sum
$$H=H_{R,-}\oplus H_{R,+},$$
where $H_{R,-}$ is the subspace generated by all eigenvectors of $A$
corresponding to eigenvalues $\lambda_{k}^{2}\leq
8(R+\lambda_{1})^{2}$, and $H_{R,+}$ is the closure of the subspace
generated by the remaining eigenvectors.  We point out that $H_{R,-}$
and $H_{R,+}$ are $A$-invariant subspaces of $H$.  Moreover, the
restriction of $A$ to $H_{R,-}$ has only a finite number of
eigenvalues, while the restriction of $A$ to $H_{R,+}$ satisfies the
coercivity condition
$$|A^{1/2}x|^{2}\geq 2(R+\lambda_{1})^{2}|x|^{2} 
\quad\quad
\forall x\in D(A^{1/2})\cap H_{R,+},$$
and even a stronger condition
\begin{equation}
	|A^{1/2}x|^{2}\geq 2(R'+\lambda_{1})^{2}|x|^{2} 
	\quad\quad
	\forall x\in D(A^{1/2})\cap H_{R,+},
	\label{hp:A-coercive-R'}
\end{equation}
for some $R'>R$ (actually we can choose $8(R'+\lambda_{1})^{2}$ to be
the smallest eigenvalue of $A$ greater than $2(R+\lambda_{1})^{2}$).

Let $u_{R,-}(t)$ and $u_{R+}(t)$ denote the components of $u(t)$ with
respect to the decomposition, let
$E(t):=|u'(t)|^{2}+|A^{1/2}u(t)|^{2}$ be the usual energy of $u(t)$,
and let $E_{R,\pm}(t)$ denote the energy of the two components.

Since the restriction of $A$ to $H_{R,-}$ has only a finite number of
eigenvalues, the component $u_{R,-}(t)$ can be regarded as a solution
to a system of finitely many ordinary differential equations.  For
this system, the time $T_{R}$ defined by (\ref{defn:TR}) coincides
with the time $t_{0}$ defined by (\ref{defn:t0-system}).  Therefore,
from Lemma~\ref{lemma:basic-u-system} (see also 
Remark~\ref{rmk:eigenspaces}) we deduce the existence of a
function $\delta:[0,T_{R}]\to[0,+\infty)$ which reduces the energy of
$u_{R,-}(t)$ at time $t=T_{R}$ by any given factor.  In particular, we
can choose this factor equal to $8e^{-2(R+\lambda_{1})T_{R}}$ and obtain that
\begin{equation}
	E_{R,-}(T_{R})\leq E_{R,-}(0)
	\cdot 8e^{-2(R+\lambda_{1})T_{R}}.
	\label{est:uR-}
\end{equation}

Now we extend $\delta(t)$ to the interval $[0,2T_{R}]$ by setting
$\delta(t):=R+\lambda_{1}$ in the second half of the interval,
namely for every $t\in(T_{R},2T_{R}]$.  Since the restriction of $A$
to $H_{R,+}$ satisfies (\ref{hp:A-coercive}), we can apply
Lemma~\ref{lemma:PDE-coercive} with $M:=R+\lambda_{1}$.  We obtain
that 
\begin{equation}
	E_{R,+}(2T_{R})\leq E_{R,+}(T_{R})
	\cdot 8e^{-2(R+\lambda_{1})T_{R}}.
	\label{est:uR+}
\end{equation}

Keeping into account that in both cases the energy is nonincreasing 
in the whole interval, from (\ref{est:uR-}) and (\ref{est:uR+}) we 
deduce that
\begin{equation}
	E(2T_{R})\leq E(0)\cdot 8e^{-2(R+\lambda_{1})T_{R}}=
	E(0)\cdot 8e^{-2\lambda_{1}T_{R}}\cdot e^{R\cdot 2T_{R}}.
	\label{est:uR}
\end{equation}

On the other hand, from (\ref{defn:TR}) it is clear that 
$2\lambda_{1}T_{R}\geq\pi$, and hence $8e^{-2\lambda_{1}T_{R}}\leq 
1$. Therefore, now (\ref{est:uR}) reads as
$$E(2T_{R})\leq E(0)\cdot e^{-R\cdot 2T_{R}},$$
which is exactly (\ref{th:u-2TR}).

If we want $\delta(t)$ with the profile of 
Figure~\ref{fig:delta-PDE}, it is enough to reduce the energy of 
$u_{R,-}(t)$ through a damping coefficient in $[0,T_{R}]$ with the 
profile shown in Figure~\ref{fig:delta-system} on the left. 

If we want a damping coefficient of class $C^{\infty}$ with all
time-derivatives vanishing at the endpoints, we need an approximation
procedure.  To this end, we first choose a bounded measurable function
$\delta(t)$ in $[0,2T_{R}]$ for which (\ref{th:u-2TR}) holds true for
all solutions with $2RT_{R}$ replaced by a larger constant.  This is
possible because in the first half of the interval we can reduce the
energy of $u_{R,-}(t)$ by any given factor, and in the second half of
the interval we can reduce the energy of $u_{R,+}(t)$ by a factor
$e^{-2(R'+\lambda_{1})T_{R}}$, where $R'>R$ is the constant that
appears in the reinforced coercivity inequality
(\ref{hp:A-coercive-R'}) (of course we need to set
$\delta(t):=R'+\lambda_{1}$ in the second half of the period, as
required by Lemma~\ref{lemma:PDE-coercive}).  At this point we can
approximate $\delta(t)$ in $L^{2}$-norm through a damping coefficient
with the required smoothness, and conclude with the aid of
Lemma~\ref{lemma:approx}.\qed

\subsubsection*{Proof of Theorem~\ref{thm:PDE-fixed}}

Let us observe that the time $2T_{R}$, with $T_{R}$ given by
(\ref{defn:TR}), coincides with $t_{0}$ as defined in
(\ref{defn:t0-PDE}).  Therefore, from Lemma~\ref{lemma:basic-PDE} we
obtain a function $\delta:[0,t_{0}]\to[0,+\infty)$ such that every
weak solution to (\ref{eqn:PDE}) satisfies
$$E(t_{0})\leq E(0)\cdot  e^{-Rt_{0}},$$
where $E(t):=|u'(t)|^{2}+|A^{1/2}u(t)|^{2}$ denotes the usual energy of 
the solution.

Now we extend $\delta(t)$ by periodicity to the whole half line
$t\geq 0$, and we obtain by an easy induction that
$$E(kt_{0})\leq E(0)\cdot e^{-kRt_{0}}
\quad\quad
\forall k\in\n.$$

Since $E(t)$ is nonincreasing, this implies that
$$E(t)\leq E(0)\cdot\exp\left(-R(t-t_{0})^{+}\right)
\quad\quad
\forall t\geq 0,$$
which is equivalent to (\ref{th:PDE-fixed}).

In analogy with the proof of Theorem~\ref{thm:ODE-fixed}, if we want a
piecewise constant damping coefficient, we can take a function with
this property from Lemma~\ref{lemma:basic-PDE}.  If we want a damping
coefficient of class $C^{\infty}$, it is enough to take from
Lemma~\ref{lemma:basic-PDE} a function of class $C^{\infty}$ whose
time-derivatives of any order vanish at the endpoints of the interval.
This condition guarantees that the periodic extension remains of class
$C^{\infty}$.\qed

\subsubsection*{Proof of Theorem~\ref{thm:PDE-top}}

Let us start by defining $\varphi(t)$. To this end, let $n_{0}$ be 
the smallest integer such that 
$\lambda_{n_{0}}^{2}>2\lambda_{1}^{2}$. For every $n\geq n_{0}$ we 
consider the positive real number
$$R_{n}:=\frac{\lambda_{n}}{\sqrt{2}}-\lambda_{1},$$
and we define $T_{n}$ as in (\ref{defn:TR}) with $R:=R_{n}$. Since 
$2(R_{n}+\lambda_{1})^{2}=\lambda_{n}^{2}$, this means that
$$T_{n}=\frac{\pi}{2}\sum_{k=1}^{n}\frac{1}{\lambda_{k}}
\quad\quad
\forall n\geq n_{0}.$$

Then we set $S_{n_{0}-1}:=0$, $U_{n_{0}-1}:=0$, and
$$S_{n}:=2\sum_{k=n_{0}}^{n}T_{k},
\hspace{4em}
U_{n}:=2\sum_{k=n_{0}}^{n}R_{k}T_{k}$$
for $n\geq n_{0}$.  It is clear that $S_{n}$ and $U_{n}$ are unbounded
increasing sequences, with
$$2T_{n+1}= 2T_{n}+\frac{\pi}{\lambda_{n+1}}\leq
2T_{n}+2T_{n_{0}}\leq S_{n}$$
for every $n\geq n_{0}+1$, and therefore
$$S_{n+1}=S_{n}+2T_{n+1}\leq 2S_{n}
\quad\quad
\forall n\geq n_{0}+1.$$

Finally, we define $\varphi:[0,+\infty)\to(0,+\infty)$ as the 
piecewise constant function such that
$$\varphi(t):=e^{-U_{n}}
\quad\quad
\mbox{if $t\in[S_{n},S_{n+1})$ for some $n\geq n_{0}-1$}.$$

It is clear that $\varphi(t)$ is nonincreasing. We claim that it 
satisfies (\ref{th:phi-decay}). Indeed for every $n\geq n_{0}-1$ it 
turns out that
\begin{equation}
	\varphi(t)e^{Rt}= e^{-U_{n}+Rt}\leq
	e^{-U_{n}+RS_{n+1}}
	\quad\quad
	\forall t\in[S_{n},S_{n+1}),
	\label{ineq:phi-eR}
\end{equation}
and for $n\geq n_{0}+1$ it turns out that
$$-U_{n}+RS_{n+1}\leq
-U_{n}+2RS_{n}=
\sum_{k=n_{0}}^{n}(4R-2R_{k})T_{k}.$$

Since $R_{k}\to+\infty$ and $T_{k}\geq T_{n_{0}}>0$, we deduce that
$-U_{n}+RS_{n+1}\to -\infty$.  Therefore, the right-hand side of
(\ref{ineq:phi-eR}) tends to 0 as $n\to+\infty$, which proves
(\ref{th:phi-decay}).

It remains to define $\delta(t)$.  To begin with, for every $n\geq
n_{0}$ we apply Lemma~\ref{lemma:basic-PDE}.  Since the time
$2T_{R_{n}}$ of Lemma~\ref{lemma:basic-PDE} coincides with $2T_{n}$ as
defined above, from the lemma we obtain a damping coefficient
$\delta_{n}:[0,2T_{n}]\to[0,+\infty)$ which reduces the energy of all
solutions to (\ref{eqn:PDE}) by a factor $e^{-2R_{n}T_{n}}$ in the
interval $[0,2T_{n}]$.

Now we glue all these functions by setting
$$\delta(t):=\delta_{n+1}(t-S_{n})
\quad\quad
\mbox{if $t\in[S_{n},S_{n+1})$ for some $n\geq n_{0}-1$.}$$

Let us consider equation (\ref{eqn:PDE}) with this choice of
$\delta(t)$, let $u(t)$ be any solution, and let
$E(t):=|u'(t)|^{2}+|A^{1/2}u(t)|^{2}$ be its usual energy.  The effect
of $\delta(t)$ in the interval $[S_{n},S_{n+1}]$ is the same as the
effect of $\delta_{n+1}(t)$ in the interval $[0,2T_{n+1}]$, and
therefore 
$$E(S_{n+1})\leq E(S_{n})\cdot e^{-2R_{n+1}T_{n+1}}
\quad\quad
\forall n\geq n_{0}-1.$$

At this point an easy induction shows that
$$E(S_{n})\leq E(0)\cdot
\exp\left(-2(R_{n_{0}}T_{n_{0}}+\ldots+R_{n}T_{n})\right)
=E(0)\cdot e^{-U_{n}}
\quad\quad
\forall n\geq n_{0}.$$

On the other hand, also for $n=n_{0}-1$ it is true that 
$E(S_{n})=E(0)\cdot e^{-U_{n}}$, for the trivial reason that 
$U_{n_{0}-1}=0$.
Since $E(t)$ is nonincreasing for $t\geq 0$, we can finally conclude 
that
$$E(t)\leq E(S_{n})\leq 
E(0)e^{-U_{n}}=E(0)\cdot\varphi(t)$$
for every $n\geq n_{0}-1$ and every $t\in[S_{n},S_{n+1}]$,
which proves (\ref{th:PDE-top}).\qed

\subsection{Further requirements on the damping coefficient}\label{sec:proofs-lip}

\subsubsection*{Proof of Proposition~\ref{prop:ODE-slow}}

Let us consider the usual energy
$E(t):=|u'(t)|^{2}+\lambda^{2}|u(t)|^{2}$.  Its time-derivative is
$E'(t)=-4\delta(t)|u'(t)|^{2}$, so that
$$E'(t)\geq -4\delta(t)E(t).$$

Integrating this differential inequality we deduce that
$$E(t)\geq E(0)\cdot\exp\left(
-4\int_{0}^{t}\delta(s)\,ds\right)
\quad\quad
\forall t\geq 0.$$

This is enough to deal with the first two cases.

For the third statement, let us assume now that $\delta(t)\geq\lambda$
for every $t\geq T$.  In this case we consider the Riccati
equation
\begin{equation}
	\varphi'(t)=\lambda^{2}-2\delta(t)\varphi(t)+\varphi^{2}(t).
	\label{eqn:ODE-riccati}
\end{equation}

Let us set $T_{*}:=\max\{T,1/\lambda\}$.  It is not difficult to check
that $\varphi(t)\equiv 0$ is a subsolution of (\ref{eqn:ODE-riccati})
for $t\geq 0$, while $\varphi(t):=\lambda-1/t$ is a supersolution for
$t\geq T_{*}$ owing to our assumptions on $\delta(t)$.  As a
consequence, the solution to (\ref{eqn:ODE-riccati}) with
$\varphi(T_{*})=0$ is defined for every $t\geq T_{*}$ and satisfies
\begin{equation}
	0\leq \varphi(t)\leq\lambda-\frac{1}{t}
	\quad\quad
	\forall t\geq T_{*}.
	\label{est:riccati-1/t}
\end{equation}

Finally, one can check that
$$u(t):=\exp\left(-\int_{T_{*}}^{t}\varphi(s)\,ds\right)
\quad\quad
\forall t\geq T_{*}$$ 
defines a solution to (\ref{eqn:ODE}), which can be easily extended to
the whole half line $t\geq 0$.  

Due to the estimate from above in (\ref{est:riccati-1/t}), this solution
satisfies 
$$|u(t)|\geq cte^{-\lambda t}
\quad\quad
\forall t\geq T_{*}$$
for a suitable positive constant $c$. Since equation (\ref{eqn:ODE}) 
is linear, $c^{-1}u(t)$ is again a solution, and satisfies 
(\ref{th:ODE-slow-lambda}).\qed
\medskip

The rest of this section is devoted to the proof of
Theorem~\ref{thm:ODE-fixed-lip}.  This requires two preliminary
results.  In the first one we show how to dampen one solution to
(\ref{eqn:ODE}) through a slowly increasing damping coefficient.

\begin{lemma}[Supercritical energy reduction of a solution]\label{lemma:lip-cut}
	Let $\lambda$ and $M$ be positive real numbers, and let 
	$\ep\in(0,\lambda)$. Let us define
	\begin{equation}
		t_{1}:=\frac{M}{2\ep}+\frac{2}{\lambda},
		\label{defn:t1}
	\end{equation}
	\begin{equation}
		\delta(t):=\lambda+\ep t
		\quad\quad
		\forall t\in[0,t_{1}].
		\label{defn:delta}
	\end{equation}
	
	Then there exists a nonzero solution $v(t)$ to (\ref{eqn:ODE}) 
	such that
	\begin{equation}
		|v'(t_{1})|^{2}+\lambda^{2}|v(t_{1})|^{2}\leq
		\left(|v'(0)|^{2}+\lambda^{2}|v(0)|^{2}\right)
		e^{-Mt_{1}}.
		\label{th:t1}
	\end{equation}
\end{lemma}

\paragraph{\textmd{\textit{Proof}}}

Let us start by proving that (\ref{defn:t1}) implies that
\begin{equation}
	(3\lambda^{2}+8\ep^{2}t_{1}^{2})
	\exp(-2\lambda t_{1}-2\ep t_{1}^{2})\leq
	2\lambda^{2}\exp(-Mt_{1}).
	\label{ineq:t1-basic}
\end{equation}

Indeed from (\ref{defn:t1}) it turns out that
$$2\lambda t_{1}\geq 4,
\hspace{4em}
2\ep t_{1}^{2}\geq Mt_{1},
\hspace{4em}
2\ep t_{1}^{2}\geq Mt_{1}+\frac{4\ep}{\lambda}t_{1}.$$

From the first two inequalities it follows that
\begin{equation}
	3\lambda^{2}\exp(-2\lambda t_{1}-2\ep t_{1}^{2})\leq
	\lambda^{2}\cdot 3e^{-4}\cdot e^{-Mt_{1}}\leq
	\lambda^{2}e^{-Mt_{1}}.
	\label{ineq:t1-basic-1}
\end{equation}

From the third inequality it follows that
$$\exp(-2\lambda t_{1}-2\ep t_{1}^{2})\leq
e^{-2\ep t_{1}^{2}}\leq e^{-4\ep t_{1}/\lambda}\cdot e^{-Mt_{1}}.$$

Therefore, from the inequality
$$t^{2}e^{-ct}\leq\frac{4}{e^{2}c^{2}}
\quad\quad
\forall t\geq 0,\ \forall c>0,$$
applied with $c:=4\ep/\lambda$ and $t:=t_{1}$, we obtain that
\begin{eqnarray}
	8\ep^{2}t_{1}^{2}\exp(-2\lambda t_{1}-2\ep t_{1}^{2}) & \leq & 
	8\ep^{2}t_{1}^{2}\exp(-4\ep t_{1}/\lambda)
	\cdot\exp(-Mt_{1})
	\nonumber  \\
	\noalign{\vspace{1ex}}
	 & \leq & 8\ep^{2}\frac{4\lambda^{2}}{16e^{2}\ep^{2}}
	 \exp(-Mt_{1})
	\nonumber  \\
	\noalign{\vspace{1ex}}
	 & \leq & \lambda^{2}\exp(-Mt_{1}).
	\label{ineq:t1-basic-2}
\end{eqnarray}

Adding (\ref{ineq:t1-basic-1}) and (\ref{ineq:t1-basic-2}) we obtain 
exactly (\ref{ineq:t1-basic}).

We are now ready to prove the existence of $v(t)$. Let us consider 
again the Riccati equation (\ref{eqn:ODE-riccati}).
When $\delta(t)$ is given by (\ref{defn:delta}), a simple computation
shows that $\varphi(t):=\lambda+2\ep t$ is a supersolution for $t\geq
0$, and the constant function $\varphi(t):=2(\lambda+2\ep t_{1})$ is a
subsolution for $t\in[0,t_{1}]$.  It follows that the solution to
(\ref{eqn:ODE-riccati}) with ``final'' condition
\begin{equation}
	\varphi(t_{1})=\lambda+2\ep t_{1}
	\label{riccati-t1}
\end{equation}
is defined for every $t\in[0,t_{1}]$ and satisfies
\begin{equation}
	\varphi(t)\geq\lambda+2\ep t
	\quad\quad
	\forall t\in[0,t_{1}].
	\label{est:phi-t}
\end{equation}

At this point a simple computation shows that
$$v(t):=\exp\left(-\int_{0}^{t}\varphi(s)\,ds\right)
\quad\quad
\forall t\in[0,t_{1}]$$
is a solution to (\ref{eqn:ODE}). We claim that this solution 
satisfies (\ref{th:t1}).

To begin with, we observe that $v(0)=1$, and from (\ref{est:phi-t}) 
it follows that
$$v(t_{1})=\exp\left(-\int_{0}^{t_{1}}\varphi(s)\,ds\right)\leq
\exp\left(-\lambda t_{1}-\ep t_{1}^{2}\right).$$

As for the time-derivative, it is given by
$$v'(t)=-\varphi(t)\exp\left(-\int_{0}^{t}\varphi(s)\,ds\right),$$
hence by (\ref{est:phi-t})
$$|v'(0)|=|\varphi(0)|\geq\lambda,$$
and by (\ref{riccati-t1}) and (\ref{est:phi-t})
$$|v'(t_{1})|=|\varphi(t_{1})|
\exp\left(-\int_{0}^{t_{1}}\varphi(s)\,ds\right)\leq
(\lambda+2\ep t_{1})\exp\left(-\lambda t_{1}-\ep t_{1}^{2}\right).$$

Recalling (\ref{ineq:t1-basic}), from all these estimates it follows 
that
\begin{eqnarray*}
	|v'(t_{1})|^{2}+\lambda^{2}|v(t_{1})|^{2} & \leq & 
	\left[(\lambda+2\ep t_{1})^{2}+\lambda^{2}\right]
	\exp\left(-2\lambda t_{1}-2\ep t_{1}^{2}\right)    \\
	\noalign{\vspace{0.5ex}}
	 & \leq & 
	 \left[3\lambda^{2}+8\ep^{2}t_{1}^{2}\right]
	 \exp\left(-2\lambda t_{1}-2\ep t_{1}^{2}\right)  \\
	 \noalign{\vspace{0.5ex}}
	 & \leq & 2\lambda^{2}\exp(-Mt_{1})  \\
	 \noalign{\vspace{0.5ex}}
	 & \leq & 
	 \left(|v'(0)|^{2}+\lambda^{2}|v(0)|^{2}\right)
	 \exp(-Mt_{1}),
\end{eqnarray*}
which proves (\ref{th:t1}).\qed
\medskip

The second preliminary result clarifies that any constant damping
coefficient below the threshold $\lambda$ allows solutions to rotate
in the phase space.

\begin{lemma}[Subcritical rotation in the phase space]\label{lemma:lip-rot}
	Let $\lambda$ be a positive real number, and let 
	$\ep\in(0,\lambda)$. Let $(\alpha,\beta)$ and $(\gamma,\delta)$ 
	be two pairs of real numbers with $\alpha^{2}+\beta^{2}\neq 0$ 
	and $\gamma^{2}+\delta^{2}\neq 0$.
	
	Then there exist positive real numbers $r$ and $t_{2}$, with
	\begin{equation}
		t_{2}\leq \frac{2\pi}{\sqrt{\ep(2\lambda-\ep)}},
		\label{th:ineq-t2}
	\end{equation}
	with the following property. The solution $u(t)$ to equation
	(\ref{eqn:ODE}) with constant dissipation 
	$\delta(t):=\lambda-\ep$ and initial data $u(0)=\alpha$ and 
	$u'(0)=\beta$ satisfies
	$$u(t_{2})=r\gamma,
	\hspace{4em}
	u'(t_{2})=r\delta.$$
\end{lemma}

\paragraph{\textmd{\textit{Proof}}}

The solution $u(t)$ can be written in the form 
$u(t):=e^{-(\lambda-\ep)t}v(t)$, where $v(t)$ is a solution to
$$v''(t)+\ep(2\lambda-\ep)v(t)=0.$$

Integrating this differential equation we see that the pair
$(v(t),v'(t))$ rotates in the phase space with period equal to the
right-hand side of (\ref{th:ineq-t2}).  During the complete rotation
the pair $(v(t),v'(t))$ turns out to be a positive multiple of any
given vector.  Since the same is true for $u(t)$, this proves the
result.\qed

\subsubsection*{Proof of Theorem~\ref{thm:ODE-fixed-lip}}

\paragraph{\textmd{\emph{Strategy}}}

\begin{figure}[htbp]
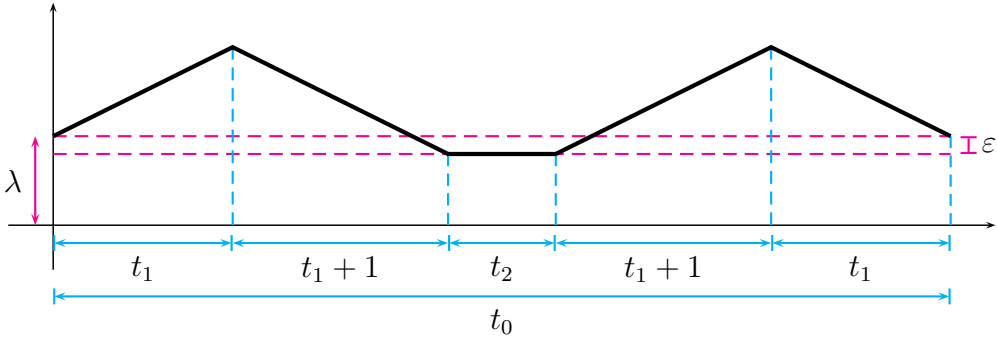

	\centering
	\psset{unit=6.5ex}
	\pspicture(-1,-1)(11,2.8)
	
	\psline[linewidth=0.7\pslinewidth]{->}(-0.5,0)(10.5,0)
	\psline[linewidth=0.7\pslinewidth]{->}(0,-0.5)(0,2.5)
	
	\psline[linestyle=dashed, linecolor=magenta](0,1)(10,1)
	\psline[linestyle=dashed, linecolor=magenta](0,0.8)(10,0.8)
	\psline[linestyle=dashed, linecolor=cyan](2,0)(2,2)
	\psline[linestyle=dashed, linecolor=cyan](4.4,0)(4.4,0.8)
	\psline[linestyle=dashed, linecolor=cyan](5.6,0)(5.6,0.8)
	\psline[linestyle=dashed, linecolor=cyan](8,0)(8,2)
	\psline[linestyle=dashed, linecolor=cyan](10,0)(10,1)
	
	\psline[linewidth=2\pslinewidth](0,1)(2,2)(4.4,0.8)(5.6,0.8)(8,2)(10,1)
	
	\pcline[offset=-0.2,linecolor=magenta]{|-|}(10,0.8)(10,1)
	\lput{U}{\uput[0](0,0){$\varepsilon$}}
	
	\pcline[offset=0.2,linecolor=magenta]{<->}(0,0)(0,1)
	\lput{U}{\uput[180](0,0){$\lambda$}}
	
	\pcline[offset=-0.2,linecolor=cyan]{|<->|}(0,0)(2,0)
	\lput{U}{\uput[270](0,0){$t_{1}$}}
	\pcline[offset=-0.2,linecolor=cyan]{<->}(2,0)(4.4,0)
	\lput{U}{\uput[270](0,0){$t_{1}+1$}}
	\pcline[offset=-0.2,linecolor=cyan]{|<->|}(4.4,0)(5.6,0)
	\lput{U}{\uput[270](0,0){$t_{2}$}}
	\pcline[offset=-0.2,linecolor=cyan]{<->}(5.6,0)(8,0)
	\lput{U}{\uput[270](0,0){$t_{1}+1$}}
	\pcline[offset=-0.2,linecolor=cyan]{|<->|}(8,0)(10,0)
	\lput{U}{\uput[270](0,0){$t_{1}$}}
	\pcline[offset=-0.8,linecolor=cyan]{|<->|}(0,0)(10,0)
	\lput{U}{\uput[270](0,0){$t_{0}$}}
	
	\endpspicture
	\caption{profile of $\delta(t)$ for the proof of Theorem~\ref{thm:ODE-fixed-lip}}
	\label{fig:delta-lip}
\end{figure}

Let us describe the strategy of the proof before entering into 
details. Let us consider the function $\delta(t)$ with the graph as 
in Figure~\ref{fig:delta-lip}, and then extended by periodicity.

The value at the endpoints is $\lambda$, the value in the horizontal 
part is $\lambda-\ep$, the slope in the oblique sections is $\pm\ep$. 
This graph depends on two parameters to be fixed, namely the length 
$t_{1}$ of the first and last oblique sections, and the length 
$t_{2}$ of the horizontal plateau. The length of each of the two 
oblique sections attached to the horizontal part is 
necessarily $t_{1}+1$. The total length is therefore
\begin{equation}
	t_{0}:=4t_{1}+t_{2}+2.
	\label{defn:t012}
\end{equation}

Now the idea is the following.  In the first ascending section,
$\delta(t)$ causes a great reduction of the energy of a special
solution $v(t)$ to (\ref{eqn:ODE}).  Let $w(t)$ be the solution with
initial data orthogonal to the initial data of $v(t)$.  In the first
two oblique sections we know that the energy of $w(t)$ is just
nonincreasing.  On the other hand, we can choose $t_{2}$ in such a way
that, at the end of the horizontal plateau, $w(t)$ has the right
``initial'' data which guarantee a reduction of its energy in the
third oblique section.

In other words, the first ascending section cuts the energy of $v(t)$,
while the horizontal plateau rotates $w(t)$ preparing it to undergo a
cut of its energy during the second ascending section.  The two
descending sections are merely junctions between the ``active'' parts
of the graph.

\paragraph{\textmd{\emph{Choice of parameters}}}

Let us set
\begin{equation}
	M:=8(\pi+\ep)R+4\ep\log 2+1,
	\label{defn:M}
\end{equation}
and let $t_{1}$ be defined by (\ref{defn:t1}). The function 
$\delta(t)$ we are considering is equal to $\lambda+\ep t$ for every 
$t\in[0,t_{1}]$. Therefore, from Lemma~\ref{lemma:lip-cut} we know 
that (\ref{eqn:ODE}) admits a solution $v(t)$ satisfying 
(\ref{th:t1}). Since the equation is linear, we can always assume that
\begin{equation}
	|v'(0)|^{2}+\lambda^{2}|v(0)|^{2}=1.
	\label{v(0)-unit}
\end{equation}

In the sequel we also need to consider $\delta(t)$ as 
extended to $t\in[-1,0]$ with the same expression $\lambda+\ep t$. 
Accordingly, we can extend $v(t)$ to the interval $[-1,0]$, and 
consider the pair $(v(-1),v'(-1))$ of its data for $t=-1$.

Let us consider now the solution $w(t)$ to (\ref{eqn:ODE}) with 
initial data
$$w(0)=\frac{v'(0)}{\lambda},
\hspace{4em}
w'(0)=-\lambda v(0).$$

We point out that the initial data of $w(t)$ satisfy
\begin{equation}
	|w'(0)|^{2}+\lambda^{2}|w(0)|^{2}=1
	\label{w(0)-unit}
\end{equation}
and they are orthogonal to the initial data of $v(t)$ in the sense 
that
\begin{equation}
	w'(0)v'(0)+\lambda^{2}w(0)v(0)=0.
	\label{v(0)-w(0)-orth}
\end{equation}

We are now ready to choose $t_{2}$. Let us consider the end of the 
second oblique section, corresponding to $t=2t_{1}+1$, and let us set
$$(\alpha,\beta):=(w(2t_{1}+1),w'(2t_{1}+1)),
\hspace{4em}
(\gamma,\delta):=(v(-1),v'(-1)).$$

From Lemma~\ref{lemma:lip-rot} we deduce the existence of a positive
time $t_{2}$ satisfying (\ref{th:ineq-t2}) such that the effect of the
constant dissipation equal to $\lambda-\ep$ on the solution with
``initial'' data $(\alpha,\beta)$ is to transform it in a multiple of
$(\gamma,\delta)$ in a time $t_{2}$. It follows that, at the end of
the horizontal plateau at time $t=2t_{1}+t_{2}+1$, the solution $w(t)$
is in the same conditions as (the extension of) the solution $v(t)$ at
time $t=-1$.  As a consequence, the energy of $w(t)$ is reduced by a
factor $e^{-Mt_{1}}$ during the third oblique section, from
$t=2t_{1}+t_{2}+1$ to $t=3t_{1}+t_{2}+2$.

\paragraph{\textmd{\emph{Estimates}}}

From (\ref{defn:t012}), (\ref{defn:t1}) and (\ref{th:ineq-t2}) it
follows that
$$t_{0}=\frac{2M}{\ep}+\frac{8}{\lambda}+2+t_{2}\leq
\frac{2M}{\ep}+\frac{8}{\lambda}+2+\frac{2\pi}{\ep},$$
so that (\ref{th:est-t0-lip}) follows from (\ref{defn:M}).

As for energy estimates, let
$E_{u}(t):=|u'(t)|^{2}+\lambda^{2}|u(t)|^{2}$ denote as usual the
energy of a solution to (\ref{eqn:ODE}).  The energy of $v(t)$ has
been reduced in the first ascending section.  Since it is always
nonincreasing, it follows that
\begin{equation}
	E_{v}(t_{0})\leq E_{v}(t_{1})\leq
	E_{v}(0)\cdot e^{-Mt_{1}}=e^{-Mt_{1}}.
	\label{est:v-t0}
\end{equation}

Similarly, the energy of $w(t)$ has been reduced in the second 
ascending section, hence
\begin{equation}
	E_{w}(t_{0})\leq E_{w}(3t_{1}+t_{2}+2)\leq
	E_{w}(2t_{1}+t_{2}+1)\cdot e^{-Mt_{1}}\leq
	E_{w}(0)\cdot e^{-Mt_{1}}= e^{-Mt_{1}}.
	\label{est:w-t0}
\end{equation}

Let us consider now any solution $u(t)$ to (\ref{eqn:ODE}). Since 
$v(t)$ and $w(t)$ are linearly independent, we can write
$$u(t)=av(t)+bw(t)$$
for suitable real constants $a$ and $b$.  Thanks to the orthonormality
relations (\ref{v(0)-unit}), (\ref{w(0)-unit}), and
(\ref{v(0)-w(0)-orth}), it turns out that
$$a^{2}+b^{2}=E_{u}(0).$$

Moreover, from (\ref{est:v-t0}) and (\ref{est:w-t0}) it follows that
\begin{eqnarray*}
	E_{u}(t_{0}) & = & 
	|av'(t_{0})+bw'(t_{0})|^{2}+\lambda^{2}|av(t_{0})+bw(t_{0})|^{2}\\
	 & \leq & 2a^{2}E_{v}(t_{0})+2b^{2}E_{w}(t_{0})  \\
	 & \leq & 2(a^{2}+b^{2})e^{-Mt_{1}}  \\
	 & = & E_{u}(0)\cdot 2e^{-Mt_{1}}.
\end{eqnarray*}

Now we claim that
\begin{equation}
	2e^{-Mt_{1}}\leq e^{-Rt_{0}},
	\label{ineq:M/R}
\end{equation}
and hence
$$E_{u}(t_{0})\leq E_{u}(0)\cdot e^{-Rt_{0}}.$$

Indeed, from (\ref{defn:M}) it follows that $M\geq 8R$, and hence
\begin{equation}
	\frac{1}{2}Mt_{1}\geq 4Rt_{1},
	\label{est:mt1-1}
\end{equation}
while from (\ref{defn:M}) and (\ref{defn:t1}) it follows that
\begin{equation}
	\frac{1}{2}Mt_{1}\geq
	\frac{M^{2}}{4\ep}\geq
	\frac{M}{4\ep}\geq
	\left(2+\frac{2\pi}{\ep}\right)R+\log 2.
	\label{est:mt1-2}
\end{equation}

Summing (\ref{est:mt1-1}) and (\ref{est:mt1-2}) we obtain that
$$Mt_{1}\geq\left(4t_{1}+2+\frac{2\pi}{\ep}\right)R+\log 2\geq
(4t_{1}+2+t_{2})R+\log 2=Rt_{0}+\log 2,$$
which is equivalent to (\ref{ineq:M/R}).

Finally, when we extend $\delta(t)$ by periodicity to the whole half 
line $t\geq 0$, an easy induction gives that
$$E_{u}(nt_{0})\leq E_{u}(0)\cdot e^{-nRt_{0}}
\quad\quad
\forall n\in\n.$$

Since $E_{u}(t)$ is nonincreasing, this implies 
(\ref{th:ODE-fixed-lip}).\qed

\setcounter{equation}{0}
\section{Applications to PDEs}\label{sec:applications}

In this section we present some simple examples of application of the 
abstract results stated in section~\ref{sec:PDE}.

Let $\Omega\subseteq\re^{d}$ be a connected bounded open set with 
smooth boundary. We consider two model examples: a dissipative wave 
equation
\begin{equation}
	u_{tt}-\Delta u+p(x)u+\delta(t)u_{t}=0
	\quad\quad
	t\geq 0,\ x\in\Omega,
	\label{eqn:wave}
\end{equation}
and a dissipative beam/plate equation
\begin{equation}
	u_{tt}+\Delta^{2} u+q(x)u+\delta(t)u_{t}=0
	\quad\quad
	t\geq 0,\ x\in\Omega,
	\label{eqn:poutre}
\end{equation}
where $p(x)$ and $q(x)$ are nonnegative bounded measurable functions. 
Just to fix ideas, we consider equation~(\ref{eqn:wave}) with homogeneous 
Dirichlet boundary conditions $u=0$, and equation~(\ref{eqn:poutre}) 
with one of the following boundary conditions: either $u=\Delta u=0$ (simply 
supported beam or plate) or $u=|\nabla u|=0$ (clamped beam or plate).

In both cases Theorem~\ref{thm:PDE-fixed} and 
Theorem~\ref{thm:PDE-top} imply that
\begin{itemize}
	\item  any exponential decay rate $e^{-Rt}$ can be realized 
	through a suitable periodic damping coefficient with period 
	$T_{R}$ that depends on $R$,

	\item with a suitable non-periodic coefficient one can achieve a
	decay rate $\varphi(t)$ faster than all exponentials.
\end{itemize}

The period $T_{R}$ and the ultra-exponential decay rate $\varphi(t)$
depend on the growth of eigenvalues, which in turn is known to depend
on the space dimension~$d$.  More precisely, when eigenvalues are
arranged in increasing order, the $n$-th eigenvalue is comparable with
the power of $n$ indicated in the fifth column of the table below (for
a proof we refer to the seminal papers~\cite{agmon,pl,weyl}).

In the third and fourth column we
state our estimates for $T_{R}$ and $\varphi(t)$.  They have to be
interpreted as asymptotic behaviors, namely when we write $\log R$ in
the third column we actually mean that $T_{R}=O(\log R)$ in that case.
In the same way, we write shortly $\exp(-t^{\alpha})$ instead of
$\exp(-O(t^{\alpha}))$.
$$
\renewcommand{\arraystretch}{1.5}
\begin{array}{|c|c||c|c||c|c|c|c|}
	\hline
	\mbox{Equation} & d & T_{R} & \varphi(t) & \lambda_{n} & T_{n} & 
	S_{n} & U_{n} \\
	\hline\hline
	\mbox{(\ref{eqn:wave})} & 1 & \log R & \exp(-t^{2}/\log t) & n & 
	\log n & n\log n & n^{2}\log n \\
	\hline
	\mbox{(\ref{eqn:wave})} & \geq 2 & R^{d-1} & \exp(-t^{2d/(2d-1)}) & 
	n^{1/d} & n^{1-1/d} & n^{2-1/d} & n^{2}  \\
	\hline
	\mbox{(\ref{eqn:poutre})} & 1 & 1 & \mbox{any} & n^2 & - & 
	- & - \\
	\hline
	\mbox{(\ref{eqn:poutre})} & 2 & \log R & \exp(-t^{2}/\log t) & n & 
	\log n & n\log n & n^{2}\log n  \\
	\hline
	\mbox{(\ref{eqn:poutre})} & \geq 3 & R^{d/2-1} & \exp(-t^{d/(d-1)}) & 
	n^{2/d} & n^{1-2/d} & n^{2-2/d} & n^{2}  \\
	\hline
\end{array}
$$

In particular, we observe that for equation (\ref{eqn:poutre}) with
$d=1$ (beam equation), the series of reciprocals of eigenvalues is
convergent.  As a consequence, the period $T_{R}$ given by
(\ref{defn:t0-PDE}) is bounded independently of $R$.  At this point,
the same argument of the proof of Theorems~\ref{thm:ODE-any}
and~\ref{thm:system-any} gives that any decay rate $\varphi(t)$ can be
achieved if we are allowed to exploit non-periodic damping
coefficients.

In all other cases, the series of reciprocals of eigenvalues is
divergent, and $T_{R}$ grows with $R$ according to
(\ref{defn:t0-PDE}).  In order to compute $\varphi(t)$, we need to
estimate the growth of the sequences $S_{n}$, $T_{n}$, $U_{n}$ defined
in the proof of Theorem~\ref{thm:PDE-top}, and recall that
$\varphi(t)$ decays in such a way that $\varphi(S_{n})=\exp(-U_{n})$.
The computations are reported in the last four columns of the table,
with the usual agreement that entries have to be intended in the sense
of the ``big~O'' notation.

The optimality of the third and fourth column of the table might 
probably deserve further future investigation.

\subsubsection*{\centering Acknowledgments}

This research project originated during a visit of the first two
authors to the Laboratoire Jacques Louis Lions of the UPMC (Paris~VI).
That stay was partially supported by the FSMP (Fondation Sciences
Math\'{e}matiques de Paris).  The project was continued later during a
visit of the third author to the University of Pisa.  That visit was
partially supported by the Gruppo Nazionale per l'Analisi
Matematica, la Probabilit\`{a} e le loro Applicazioni (GNAMPA) of the
Istituto Nazionale di Alta Matematica (INdAM).

\label{NumeroPagine}

\end{document}